\documentclass[11pt]{amsart}
\usepackage[T1]{fontenc}
\usepackage{amsmath,amssymb,mathtools,booktabs,tikz-cd,mathrsfs, comment}
\usepackage[margin=1.1in]{geometry}
\usepackage[colorlinks=true,linkcolor=blue,citecolor=blue,urlcolor=blue]{hyperref}
\newtheorem{theorem}{Theorem}[section]
\newtheorem{proposition}[theorem]{Proposition}
\newtheorem{lemma}[theorem]{Lemma}
\newtheorem{corollary}[theorem]{Corollary}
\theoremstyle{definition}

\newtheorem{conjecture}[theorem]{Conjecture}

\theoremstyle{remark}
\newtheorem{remark}[theorem]{Remark}
\newtheorem{example}[theorem]{Example}

\newcommand{\C}{\mathbb C}
\newcommand{\Z}{\mathbb Z}
\newcommand{\g}{\mathfrak g}
\newcommand{\n}{\mathfrak n_+}
\newcommand{\h}{\mathfrak h}
\newcommand{\vac}{\mathbf1}
\newcommand{\J}{J_\infty}
\newcommand{\spalg}{\mathfrak{sp}}
\newcommand{\slalg}{\mathfrak{sl}}
\newcommand{\soalg}{\mathfrak{so}}
\newcommand{\normal}[1]{\mathopen{:}#1\mathclose{:}}
\DeclareMathOperator{\gr}{gr}
\DeclareMathOperator{\Spec}{Spec}
\DeclareMathOperator{\ad}{ad}
\DeclareMathOperator{\rank}{rank}
\DeclareMathOperator{\ch}{ch}

\DeclareMathOperator{\wt}{wt}

\newcommand{\B}{\mathcal B}
\newcommand{\SingSupp}{\operatorname{SS}}
\newcommand{\HS}{\operatorname{Hilb}}
\newcommand{\red}{\mathrm{red}}
\newcommand{\slthree}{\mathfrak{sl}_3}
\title[Principal subspaces at admissible levels]{Principal subspaces of affine vertex algebras at admissible levels}
\author{Hao Li}
\email{hao-li@oist.jp}
\author{Antun Milas}
\email{amilas@albany.edu}
\date{}
\subjclass[2020]{Primary 17B69; Secondary 14E18, 17B67, 05A17}
\keywords{Principal subspace, admissible affine vertex algebra, nilpotent orbit, Li filtration, arc algebra, classical freeness}

\begin{document}
\begin{abstract}
We classify generic principal subspaces of admissible affine
vertex algebras for all simple Lie algebras. The nongeneric cases occur
precisely at denominator two in types $B,C,F_4$ and at denominator three
in type $G_2$. The proof uses geometry of nilpotent orbits combined with positive root
null vectors. Surprisingly, null vectors needed for our classification are explicitly known and already 
studied in the literature.

We then study the exceptional principal subspaces at the
level of their $C_2$-algebras and arc algebras. For $\mathfrak{sp}_4$ at
level $-\frac12$ and $G_2$ at level $-\frac53$, we prove classical freeness
and obtain product character formulas.  In type $C_n$ at level
$-\frac12$, classical freeness fails for $n\ge4$. For $n=3$, we combine the functional model method of Feigin and Stoyanovsky, and a flat degeneration of the explicit realization of $L_{-3/2}(\mathfrak{sl}_3)$ to prove a product character formula, while classical freeness in this case remains conjectural. The principal subspace of $L_{-5/2}(F_4)$ is also conjecturally 
classically free.
\end{abstract}
\maketitle

\section{Introduction}
Principal subspaces connect affine Lie algebra representations, partition
identities, and commutative algebra. Introduced by Feigin and Stoyanovsky
\cite{FS}; see also \cite{CLM0, CLM,LiMilas}, they are generated from a highest-weight vector by positive-root
generators. Integrability imposes relations among these generators at positive
integral level. At nonintegral admissible level, the ambient affine vertex
algebra is again a proper quotient of its universal counterpart, but its
defining ideal need not meet the universal principal subspace.

Fix a triangular decomposition of a simple Lie algebra $\g$, with positive
nilpotent subalgebra $\n$. We write
\begin{equation}\label{eq:intro-principal}
 W(k\Lambda_0)=U(\n\otimes t^{-1}\C[t^{-1}])\vac\subset L_k(\g),
 \qquad N(k\Lambda_0)\subset V^k(\g)
\end{equation}
for the simple and universal principal subspaces, respectively; the Lie
algebra is specified by the context. Since the invariant form vanishes on
$\n$, $N(k\Lambda_0)\cong V^0(\n)$ is independent of $k$. We call $W(k\Lambda_0)$
\emph{generic} when the canonical map $N(k\Lambda_0)\to W(k\Lambda_0)$
is an isomorphism. 

Our first main result is classification of such cases when the level is admissible.
\begin{theorem}\label{thm:intro-classification}\label{thm:classification}
Let $k=-h^\vee+p/q$ be a nonintegrable admissible vacuum level in lowest
terms, allowing both principal and coprincipal admissibility. Then
$W(k\Lambda_0)$ is nongeneric precisely in the following cases:
\[
\begin{array}{c|c|l}
\g&q&k\\\hline
B_\ell\ (\ell\ge2)&2&-\ell+\frac32+j\\
C_n\ (n\ge3)&2&-\frac12+j\\
F_4&2&-\frac52+j\\
G_2&3&-\frac53+j\ \text{or}\ -\frac43+j
\end{array}
\qquad (j\in\Z_{\ge0}).
\]
In particular, all nonintegrable admissible principal subspaces in types
$A,D,E$ are generic. In every generic case,
\begin{equation}\label{eq:intro-generic-char}
 \ch W(k\Lambda_0)(q)=\prod_{r\ge1}\frac{1}{(1-q^r)^{|\Delta_+|}}.
\end{equation}
\end{theorem}

The exceptional families are coprincipal admissible in the standard
terminology; restricting to principal admissible levels would exclude
them. Our classification permits every admissible numerator, rather than
only boundary levels \cite{KW}.

The geometric input needed in the paper is not really new because there are numerous results already applicable
to affine vertex algebras at admissible level. We basically applied such results to principal subspaces.
The classification combines Arakawa's description of associated varieties with Panyushev's generic stabilizers and known null vectors.

 Arakawa's identification is
$X_{L_k(\g)}=\overline{\mathbb O_q}$ \cite{ArAssoc,AvEM}.
Throughout the paper, $\pi:\g^*\to\n^*$ denotes restriction of linear
functionals; its restriction to a nilpotent orbit is the moment map for
$N_+$. Dominance of this restriction implies genericity of the principal
subspace, using Li's filtration \cite{Li} and the reduced
singular-support theorem for quasi-lisse vertex algebras \cite{AMArc}.

The necessary generic stabilizers are computed directly or obtained from
Panyushev \cite{Pan}. Conversely, positive-root null vectors in types
$B,C,F_4$ \cite{PerseB,AMJoseph,PerseF} and $G_2$ \cite{AL} establish
nongenericity. It is remarkable that in all these cases null vectors are 
computed a while ago. It is easy to see that non-genericity propagates under
nonnegative integral shifts of the level. This proves Theorem \ref{thm:intro-classification}. 

We next study the canonical arc map
\begin{equation}\label{eq:intro-arcmap}
 \psi_V:\J(R_V)\twoheadrightarrow\gr^F V,
 \qquad R_V=V/C_2(V).
\end{equation}
Classical freeness means that this map is an isomorphism: the relations in
the Li associated graded algebra are generated by the $C_2$-relations and their
derivatives. Generic principal subspaces are classically free by PBW.
Among the exceptional cases, we prove
\[
\begin{array}{c|c|c}
\g, k&R_{W(k\Lambda_0)}&\ch W(k\Lambda_0)\\\hline
\mathfrak{sp}_4,\ -\frac12&
\C[x,y,z,u]/(xy-z^2)&\dfrac1{(1-q)(q;q)_\infty^3}\\[5pt]
G_2,\ -\frac53&
\C[a,b,c,d,e,f]/(\frac13d^2-ce+af)&
\dfrac1{(1-q)(q;q)_\infty^5}
\end{array}
\]
and classical freeness in both cases. The proofs use the geometry of
quadric jet schemes and injectivity of the corresponding arc maps. 
For $\mathfrak{sp}_4$ the result extends to all nongeneric levels. For every
$\ell\ge 0$ we prove that $W((\ell-\tfrac12)\Lambda_0)$ is classically free, with
\[
R_{W}\cong \mathbb C[x,y,z,u]/\bigl((xy-z^2)^{\ell+1}\bigr),\qquad
\operatorname{ch}W((\ell-\tfrac12)\Lambda_0)=\frac{1}{(1-q)(q;q)_\infty^3}
\prod_{\substack{m\ge 1\\ m\not\equiv 0,\pm1 \ (\mathrm{mod}\ 2\ell+3)}}\frac{1}{1-q^m},
\]
so the character is governed by the Andrews--Gordon identities
(Theorem~\ref{thm:sp4-half-integral}). For $G_2$ at the second initial level $-\tfrac43$,
the square of the quadratic relation gives a surjection onto $R_W$. We conjecture
that it is an isomorphism and that $W(-\tfrac43\Lambda_0)$ is classically free, which
would give a character of Rogers--Ramanujan type (Conjecture~\ref{conj:rogers-ramanujian}).

The symplectic-boson realization \cite{FF} behaves differently in rank
three. The intrinsic $C_2$-algebra has a nonzero nilpotent class $h$ invisible
in the ambient free-field quotient $R_{\mathcal{B}_3}$. 
We show that eight quadratic
and three cubic relations present it:
\[
R_W\cong\widetilde{T}[d],\qquad
\widetilde{T}=\mathbb{C}[x,y,z,p,q,r,a,b]/\bigl(I_2+(ah,bh,rh)\bigr).
\]
Any further relation would be of the form $hd^m=0$, and this is excluded by a symbol computation
in the free-field realization. We then prove
\begin{equation}
	\operatorname{ch}W(-\tfrac12\Lambda_0)=\frac{(-q;q)^8_\infty}{(q;q)_\infty}.
\end{equation}
Removing the generator $D_{12}$ leaves an eight-generator subalgebra $\mathcal{A}$ with
$\operatorname{ch}W=\operatorname{ch}\mathcal{A}/(q;q)_\infty$, so (4) amounts to
$\operatorname{ch}\mathcal{A}=(-q;q)^8_\infty$. For the lower bound, let
$\widehat{S}\subset J_\infty(R_{B_3})$ be the differential subalgebra generated by the symbols of
the eight generators and of the state representing $h$. Then
$\operatorname{ch}\mathcal{A}\ge\operatorname{Hilb}\widehat{S}$. A functional model method originated from Feigin--Stoyanovsky \cite{FS} reduces $\operatorname{Hilb}\widehat{S}$ to the rank-two arc
algebra of Section~5 and $q$-binomial coefficients, and the resulting sum is evaluated by Rogers' $q$-Mehler formula. For the upper bound, motivated by the
deformation methods of Stoyanovsky \cite{Stoy} further developed by Creutzig and Garner \cite{CG}, the currents of $\mathcal{A}$ span a contraction of
$\mathfrak{sl}_3$. Degenerating the free-field realization of $L_{-3/2}(\mathfrak{sl}_3)$
\cite{Ad16,ACG,BMR} gives a surjection onto $\mathcal{A}$ from the special fibre of a flat degeneration of
$L_{-3/2}(\mathfrak{sl}_3)$, whose character is $(-q;q)^8_\infty$ by Kac--Wakimoto \cite{KW}. The
two bounds agree, which proves (4). Consequently, $W$ is classically free if and only if
\[
\operatorname{Hilb}J_\infty(\widetilde{T})\le(-q;q)^8_\infty
\]
coefficientwise, and in that case $L_{-3/2}(\mathfrak{sl}_3)$ is classically free as well.

\color{black}

\color{black}

For $\g=\mathfrak{sp}_{2n}$ and $n\ge4$, a computer-assisted rank-four
calculation and a terminal-Levi-face argument prove failure of classical
freeness at level $-\frac12$. The analogous question for
$W(-\frac52\Lambda_0)$ in type $F_4$ is left open, with an explicit
candidate presentation suggested by the projected minimal orbit.
We also conjecture, for all admissible levels, that the reduced associated variety
of the principal subspace is $\overline{\pi(\mathbb O_q)}$
(Conjecture~\ref{conj:principal-associated-variety}). This geometric conjecture allows nilpotents in the
intrinsic $C_2$-algebra and is independent of classical freeness. It holds in all
generic cases, and our case studies confirm it for $G_2$ at level $-\tfrac53$ and for
$\mathfrak{sp}_6$ at level $-\tfrac12$, where $R_W$ is not reduced.

Section~\ref{sec:prelim} fixes the common notation and filtration tools.
Section~\ref{sec:classification} proves the classification, and
Section~\ref{subsec:associated-variety-conjecture} formulates the
associated-variety conjecture. The single-section case studies for
$\mathfrak{sp}_4$, $G_2$, and $\mathfrak{sp}_6$ follow in
Sections~\ref{sec:main}, \ref{sec:g2-classical}, and
\ref{sec:sp6}. Section~\ref{sec:higher-rank} treats the higher-rank
obstruction, and Section~\ref{subsec:F4-conjecture} discusses $F_4$.


(Use of AI) This paper originates from our 2021 preprint in which several results of this paper were written. Our initial draft had an incomplete proof/conjecture of classical freeness of the principal subspace $\mathfrak{sp}(2n)_{-1/2}$, using the free field construction. We recently ran an experiment with $\mathfrak{sp}(8)_{-1/2}$ with ChatGPT and found a mistake in our argument.  
In Section 7, the idea of the proof was suggested by ChatGPT; the proof based on this method was written by the authors. The OPE computations for the degeneration of the explicit realization of $L_{-3/2}(\mathfrak{sl}_3)$ were carried out by ChatGPT and verified by the authors.  After drafting our results, we edited and proofread the whole paper with ChatGPT. 
\section*{Acknowledgements}
H.~Li thanks the Okinawa Institute of Science and Technology Graduate University (OIST) for
providing access to ChatGPT~Edu.
A.~Milas received partial support from the Simons Foundation (Travel Support for
Mathematicians, grant number 709563) and the U.S.\ National Science Foundation (grant number
DMS-2101844). The work of A.~Milas was also funded by the European Union -- NextGenerationEU
through the Juraj Dobrila University of Pula institutional research project
IIP\_UNIPU\_010159.

\section{Notation, filtrations, and arc algebras}\label{sec:prelim}
\subsection{Principal subspaces and grading}
Fix $\g=\mathfrak n_-\oplus\h\oplus\n$ and normalize the invariant
form so that long roots have squared length two. Put $a(r)=a\otimes t^r$;
then the affine bracket at level $k$ is
\[
 [a(r),b(s)]=[a,b](r+s)+r k(a,b)\delta_{r+s,0}.
\]
Write $\widehat\n=\n\otimes\C[t,t^{-1}]$ and
$\mathfrak u=\n\otimes t^{-1}\C[t^{-1}]$.
The principal subspaces in \eqref{eq:intro-principal} are vertex
subalgebras, strongly generated by the generators of a basis of $\n$.
We use the canonical quotient
\begin{equation}\label{eq:principal-map}
 \rho_k:N(k\Lambda_0)\twoheadrightarrow W(k\Lambda_0).
\end{equation}
The notation $N(k\Lambda_0)$ records its ambient vacuum module;
its abstract vertex algebra is $V^0(\n)$ for every $k$.
In a fixed case we abbreviate these algebras to $N$ and $W$.
The kernel of \eqref{eq:principal-map} is a vertex ideal, whereas the
kernel of $U(\widehat\n)\to W$, $u\mapsto u\vac$, is a left ideal.

The inherited affine degree of
$a_1(-m_1)\cdots a_r(-m_r)\vac$ is $\sum_i m_i$, with $m_i\ge1$.
We use this grading even when the ambient conformal vector does not
belong to $W$. Characters are unshifted. In
Sections~\ref{sec:classification}--\ref{subsec:associated-variety-conjecture},
$q$ is the admissible denominator and $z$ is used for characters;
in the subsequent case studies, $q$ is the character variable.
We write $(a;q)_\infty=\prod_{r\ge0}(1-aq^r)$.

\subsection{Li's filtration and classical freeness}
For a vertex algebra $V$, Li's decreasing filtration is
\[
 F^pV=\operatorname{span}\left\{
 a^1_{(-r_1-1)}\cdots a^s_{(-r_s-1)}b:
 r_i\ge0,\ \sum_i r_i\ge p\right\}.
\]
Set $R_V=V/F^1V=V/C_2(V)$ and
$\gr^F V=\bigoplus_{p\ge0}F^pV/F^{p+1}V$.
For the affine Lie algebras above, $F^p$ is spanned by negative-current
monomials with $\sum_i(m_i-1)\ge p$.
The associated graded is a commutative differential algebra, and
$\SingSupp(V)=\Spec(\gr^F V)$ is the singular-support scheme.
We write $X_V=\operatorname{Specm}R_V$ for the reduced associated variety introduced in \cite{Ar}.

For $R=\C[x_1,\ldots,x_s]/I$, set
$x_a(t)=\sum_{r\ge0}x_a^{(r)}t^r$ and define
\begin{equation}\label{eq:jet-definition}
 J_m(R)=\frac{\C[x_a^{(r)}:1\le a\le s,\ 0\le r\le m]}
 {([t^r]f(x(t)):f\in I,\ 0\le r\le m)},\qquad
 \J(R)=\varinjlim_mJ_m(R).
\end{equation}
The translation derivation is $Tx_a^{(r)}=(r+1)x_a^{(r+1)}$.
If $x_a$ has weight $d_a$, then $x_a^{(r)}$ has weight $d_a+r$.
Finite jets impose only the displayed coefficient equations, not all
iterates of a derivation truncated at order $m$.
There is a canonical surjection
\begin{equation}\label{eq:canonical}
 \psi_V:\J(R_V)\twoheadrightarrow\gr^F V,
 \qquad \bar a^{(r)}\longmapsto\sigma_r(T^ra/r!).
\end{equation}
Here $\sigma_r$ denotes the class in $F^rV/F^{r+1}V$.
We call $V$ \emph{classically free} when this map is an isomorphism.
PBW gives
\begin{equation}\label{eq:universal-pbw}
 R_{N(k\Lambda_0)}\cong S(\n),\qquad
 \gr^F N(k\Lambda_0)\cong\J(S(\n)).
\end{equation}
All filtrations used here are finite in each affine degree, so passage
to the associated graded preserves the character.

\subsection{An injectivity criterion}\label{sec:criterion}
Let $j:A\hookrightarrow B$ be an embedding of vertex algebras, and let
$\bar j:R_A\to R_B$ be the induced map. Functoriality gives
\begin{equation}\label{eq:diagram}
\begin{tikzcd}[column sep=large,row sep=large]
\J(R_A)\arrow[r,two heads,"\psi_A"]\arrow[d,"\J(\bar j)"']&
\gr^F A\arrow[d,"\gr^Fj"]\\
\J(R_B)\arrow[r,two heads,"\psi_B"]&\gr^F B.
\end{tikzcd}
\end{equation}
\begin{proposition}\label{thm:criterion}
If $B$ is classically free and $\J(\bar j)$ is injective, then $A$ is
classically free and $\gr^Fj$ is injective.
\end{proposition}
\begin{proof}
The composite $\psi_B\J(\bar j)=(\gr^Fj)\psi_A$ is injective.
Since $\psi_A$ is surjective, both assertions follow.
\end{proof}
\begin{remark}\label{rem:strictness} \cite{AL1}
Neither $\bar j$ nor $\gr^Fj$ is assumed injective merely because $j$ is.
Under the proposition's hypotheses, degree-wise finiteness of the Li
filtration also gives $F^pA=A\cap F^pB$ by induction on filtration degree.
An injection of coordinate rings alone is insufficient: under
$\C[x,y]/(y^2-x^3)\hookrightarrow\C[t]$, $x\mapsto t^2$, $y\mapsto t^3$,
the nonzero arc element $2x\,Ty-3y\,Tx$ maps to zero. Its nonvanishing
follows by assigning weights two and three to all derivatives of $x$ and
$y$: it has weight five, while the differential ideal is generated in
weight six. 
\end{remark}


\section{Generic principal subspaces at admissible levels}
\label{sec:classification}

\subsection{Admissible levels}
We prove Theorem~\ref{thm:classification} using the notation of
Section~\ref{sec:prelim}.
Let $r^\vee$ be the lacing number. The admissible vacuum levels are
\begin{equation}\label{eq:admissibility}
 k=-h^\vee+\frac pq,\quad p,q\in\Z_{>0},\quad (p,q)=1,
 \qquad
 \begin{cases}
 p\ge h^\vee,&(q,r^\vee)=1,\\
 p\ge h,&(q,r^\vee)=r^\vee.
 \end{cases}
\end{equation}
For the simple types, these two cases exhaust the possibilities.
They are called principal and coprincipal admissible, respectively
(with no distinction in simply-laced type). In particular, restricting
to principal admissible levels alone would exclude the denominator-two
and denominator-three exceptions considered below.

\begin{remark}
The restriction to nonintegrable levels is of course necessary. At an integrable
level $k\in\Z_{\ge0}$, the relation
$e_\theta(-1)^{k+1}\vac=0$ makes the principal subspace nongeneric.
\end{remark}

\subsection{A moment-map criterion}
Let $G$ be the connected adjoint group of $\mathfrak g$, and identify
$\mathfrak g\simeq\mathfrak g^*$ using a nondegenerate invariant symmetric
bilinear form $(\,\cdot\,,\,\cdot\,)$. 
By
\cite{ArAssoc,AvEM}, the reduced associated variety is
$X_{L_k(\g)}=\overline{\mathbb O_q}$, where
\begin{equation}\label{eq:orbit-description}
 \overline{\mathbb O_q}=
 \begin{cases}
 \{x\in\g:(\ad x)^{2q}=0\},&(q,r^\vee)=1,\\
 \{x\in\g:\varrho_{\theta_s}(x)^{2q/r^\vee}=0\},
 &(q,r^\vee)=r^\vee.
 \end{cases}
\end{equation}
Here $\varrho_{\theta_s}$ is the irreducible representation of highest
short root. 

For a nilpotent orbit $\mathbb O$, set
\[
 d(\mathbb O)
 =\min_{x\in\mathbb O}\dim(\mathfrak g^x\cap\mathfrak n_+),
 \qquad
 \mathfrak g^x=\{a\in\mathfrak g:[a,x]=0\}.
\]
Let $\pi:\mathfrak g^*\to\mathfrak n_+^*$ be restriction of functionals.
Since $T_x\mathbb O=[\mathfrak g,x]$, its differential satisfies
\[
 \bigl(d(\pi|_{\mathbb O})_x([a,x])\bigr)(u)
 =([a,x],u)=(a,[x,u]),
 \qquad a\in\mathfrak g,\quad u\in\mathfrak n_+.
\]
Nondegeneracy therefore gives
\[
 \ker\bigl(d(\pi|_{\mathbb O})_x\bigr)^*
 =\mathfrak g^x\cap\mathfrak n_+,
\]
and hence
\begin{equation}\label{eq:moment-rank}
 \operatorname{rank}d(\pi|_{\mathbb O})_x
 =\dim\mathfrak n_+-\dim(\mathfrak g^x\cap\mathfrak n_+)
 =\dim(N_+\cdot x).
\end{equation}

This is also the moment-map rank formula. Indeed, with the convention
\[
 \omega_x([a,x],[b,x])=(x,[a,b])
\]
for the Kirillov--Kostant form, the function
$\pi^u(x)=\pi(x)(u)=(x,u)$ satisfies
\[
 d\pi^u_x([a,x])
 =(x,[u,a])
 =\omega_x([u,x],[a,x]).
\]
Thus $\pi|_{\mathbb O}$ is the moment map for the $N_+$-action,
and the symplectic form identifies the transpose of its differential
with the infinitesimal action $u\mapsto[u,x]$.

Over $\mathbb C$, the dimension of the image closure equals the generic,
equivalently maximal, differential rank. Consequently,
\begin{equation}\label{eq:projection-dimension}
 \dim\overline{\pi(\mathbb O)}
 =\dim\mathfrak n_+-d(\mathbb O).
\end{equation}
Here $d(\mathbb O)$ is the stabilizer dimension on a dense open subset
of $\mathbb O$; no assertion about conjugacy of these stabilizers is
needed. Since a zero-dimensional algebraic subgroup of a complex
unipotent group is trivial, $\pi|_{\mathbb O}$ is dominant if and only if
the $N_+$-action on $\mathbb O$ is generically free, that is,
\[
 (N_+)_x
 =\{n\in N_+:\operatorname{Ad}(n)x=x\}
 =\{1\}
\]
for $x$ in a dense open subset of $\mathbb O$, and thus the restriction is
dominant exactly when $d(\mathbb O)=0$.


 The principal subspace $W(k\Lambda_0)$ is generic when the
positive-root modes $u(-n)$, $u\in\mathfrak n_+$, $n\ge1$, satisfy no relations in
$L_k(\mathfrak g)$, that is, when $\rho_k$ is injective. In the $C_2$-algebra only the modes
$u(-1)$ survive, and their symbols restrict to the linear functions $\ell_u(x)=(x,u)$ on the
associated variety $\overline{\mathbb O_q}$. The classical shadow of genericity is therefore that the
$\ell_u$ satisfy no further relations on $\overline{\mathbb O_q}$; equivalently, that
$\pi|_{\overline{\mathbb O_q}}$ is dominant, which by \eqref{eq:moment-rank} means $d_q=0$. 
To control all
modes $u(-n)$, we have pass to arc spaces and use additional information. 
\color{black}

\begin{lemma}\label{lem:semicontinuity}
If $\mathbb O'\subset\overline{\mathbb O}$, then
$d(\mathbb O)\le d(\mathbb O')$.
\end{lemma}
\begin{proof}
The function $x\mapsto\dim\ker(\n\to\g,\ u\mapsto[u,x])$
is upper semicontinuous. An open set on which its value is at most
$d(\mathbb O')$ meets the dense orbit $\mathbb O$.
\end{proof}

Next result is certainly known to experts but we include proof for completeness.
\begin{lemma}\label{lem:arc-dominance}
Let $f\colon X\to Y$ be a dominant morphism of irreducible affine varieties over $\mathbb C$,
and assume that every jet scheme $J_mY$ is reduced and irreducible (for example $Y=\mathbb A^s$).
Then the pullback $\mathbb C[J_\infty Y]\to\mathbb C[J_\infty X]_{\mathrm{red}}$ is injective.
\end{lemma}
\begin{proof}
By generic smoothness there are nonempty open subsets $U\subset X_{\mathrm{sm}}$ and $V\subset Y$
such that $f|_U\colon U\to V$ is smooth and surjective. Let $\gamma$ be an $m$-jet of $Y$ centred
in $V$. Choose a point of $U$ over $\gamma(0)$; by formal smoothness, $\gamma$ lifts to an $m$-jet
of $U$, and this jet extends to an arc because $U$ is smooth. Every element of
$\mathbb C[J_\infty Y]$ is the image of some $g\in\mathbb C[J_mY]$. If its pullback is nilpotent,
then $g$ vanishes on all $m$-jets centred in $V$. These form a nonempty open subset of the
irreducible scheme $J_mY$, hence a dense one. Since $J_mY$ is reduced, $g=0$.
\end{proof}

\color{black}

\begin{proposition}\label{prop:genericity}
If $k$ is admissible and $d_q=0$, then $W(k\Lambda_0)$ is generic.
\end{proposition}
\begin{proof}
Put $V=L_k(\g)$ and $X=\overline{\mathbb O_q}$.
Since $X$ has finitely many symplectic leaves, $V$ is quasi-lisse.
The canonical reduced singular-support identification
\cite[Theorem~9.2]{AMArc} is
\[
 \eta:\C[J_\infty X]_{\mathrm{red}}
 \xrightarrow{\sim}(\gr^F V)_{\mathrm{red}}.
\]
Consider the induced composite
\begin{equation}\label{eq:graded-composite}
 \J(S(\n))=\gr^F N(k\Lambda_0)
 \longrightarrow\gr^F V
 \longrightarrow(\gr^F V)_{\mathrm{red}}.
\end{equation}
For $u\in\n$, put $\ell_u(x)=(x,u)$. The pullback $\pi^*$ sends
$u$ to $\ell_u$, and its arc map sends $u^{(r)}$ to $\ell_u^{(r)}$.
By construction, $\eta$ sends $\ell_u^{(r)}$ to the reduced Li symbol of
$u(-r-1)\vac=T^r(u(-1)\vac)/r!$. This is also the image of $u^{(r)}$
under \eqref{eq:graded-composite}. Consequently the diagram
\[
\begin{tikzcd}[column sep=large,row sep=large]
 \J(S(\n))\arrow[r,"\J(\pi^*)"]
 \arrow[dr,"\eqref{eq:graded-composite}"']&
 \C[J_\infty X]_{\mathrm{red}}\arrow[d,"\eta","\sim"']\\
 &(\gr^F V)_{\mathrm{red}}
\end{tikzcd}
\]
commutes on every jet generator. If $d_q=0$, the top arrow is injective
by \eqref{eq:moment-rank} and Lemma~\ref{lem:arc-dominance}.
Thus \eqref{eq:graded-composite} is injective. A nonzero homogeneous
element of $\ker\rho_k$ would have a nonzero leading Li symbol in its
kernel, since the filtration on $N(k\Lambda_0)$ is finite in each degree.
Hence $\ker\rho_k=0$.
\end{proof}


\begin{lemma}\label{lem:propagation}
Suppose that $W(k\Lambda_0)$ is nongeneric and let $\ell\in\Z_{\ge0}$.
Assume that the simple affine algebras at the levels in question are
defined as graded vacuum quotients. Then $W((k+\ell)\Lambda_0)$ is
nongeneric. In particular, this applies to admissible $k$ and its
nonnegative integral translates.
\end{lemma}
\begin{proof}
Since $W(k\Lambda_0)$ is nongeneric there exists $0 \neq P \in U(\widehat{\frak n_+}_{< 0})$ such that $P \vac=0$. Let
$
\mathfrak u=\mathfrak n_+\otimes t^{-1}\mathbb C[t^{-1}],
$
and let \(\mathfrak a\subset\mathfrak u\) be the finite-dimensional nilpotent Lie algebra generated by the root modes occurring in \(P\). Put
\[
d=\dim\mathfrak a,\qquad m=\deg(P),\qquad
E=U(\mathfrak a)\mathbf1\subset L_\ell(\mathfrak g).
\]
Integrability implies that $E$ is finite-dimensional.
Denote by $U(\mathfrak a)_{\le m}$ PBW  filtration of $U(\mathfrak a)$
Then we have to show that there exists $r$ large enough and 
$0\ne Q\in U(\mathfrak a)_{\le r}
\text{ with }
Q(\mathbf1\otimes\mathbf1)=0,
$
where \(Q\) acts diagonally. Consequently,
\[
Q\mathbf1=0\quad\text{in } \quad W_{k+\ell}(\frak g) \subset L_{k+\ell}(\mathfrak g).
\]

To prove this, write  \(M=U(\frak a) \mathbf1\subset L_k(\mathfrak g)\). Since \(UP\) annihilates the vacuum, using combinatorial formula for counting 
basis PBW vectors gives, for \(r\ge m\),
\[
\dim U_{\le r}\mathbf1
\le
\binom{r+d}{d}-\binom{r-m+d}{d}.
\]
Indeed, right multiplication by \(P\) embeds \(U_{\le r-m}\) into \(U_{\le r}\), because \(U\) is a domain. Under the diagonal action,
\[
U_{\le r}(\mathbf1\otimes\mathbf1)
\subseteq (U_{\le r}\mathbf1)\otimes E.
\]
Therefore its dimension is at most
\[
{\rm dim}(E) \cdot \left[\binom{r+d}{d}-\binom{r-m+d}{d}\right]
\le
\frac{{\rm dim}(E) \cdot md}{r+d}\binom{r+d}{d}.
\]
Taking \(r={\rm dim}(E)md\) makes this strictly smaller than \(\dim U_{\le r}\). Hence a nonzero \(Q\) exists. \end{proof}

\begin{remark} In the case when the Lie algebra $\frak a$ is commutative, we can obtain a stronger result, namely $P^r {\bf 1}=0$ for some $r$. This explains why some null vectors propagate as higher powers of null vectors.
\end{remark}

\subsection{Type \texorpdfstring{$A$}{A}}
\begin{proposition}\label{prop:typeA}
For $\g=\slalg_N$, $N\ge2$, one has $d_q=0$ for every $q\ge2$.
\end{proposition}
\begin{proof}
Put $m=\lfloor N/2\rfloor$ and
$x_0=\sum_{i=1}^m E_{N+1-i,i}$. Then $x_0^2=0$. For $N=2m$,
write $J$ for the antidiagonal identity matrix of size $m$ and
\[
 x_0=\begin{pmatrix}0&0\\J&0\end{pmatrix},\qquad
 u=\begin{pmatrix}A&B\\0&D\end{pmatrix}\in\n.
\]
The equation $[u,x_0]=0$ implies $B=0$ and $D=JAJ$.
Since $A,D$ are strictly upper triangular and $JAJ$ is strictly
lower triangular, $u=0$. For $N=2m+1$, insert a middle block of size
one:
\[
 x_0=\begin{pmatrix}0&0&0\\0&0&0\\J&0&0\end{pmatrix},\qquad
 u=\begin{pmatrix}A&a&B\\0&0&b\\0&0&D\end{pmatrix}.
\]
Commutation also forces $a=b=0$, and again $u=0$.
Finally, $x_0^2=0$ implies $(\ad x_0)^3y=x_0^3 y-3x_0^2 y x_0+3x_0 y x_0^2-y x_0^3 =0$ for every $y$.
Thus $x_0\in\overline{\mathbb O_q}$ for every $q\ge2$ by
\eqref{eq:orbit-description}, and Lemma~\ref{lem:semicontinuity}
finishes the proof.
\end{proof}

\subsection{Type \texorpdfstring{$C$}{C}}
\begin{proposition}\label{prop:typeC}
For $\g=\spalg_{2n}$, $n\ge2$, one has
\[
 d_2=(n-1)^2,\qquad d_q=0\quad(q\ge3).
\]
The principal subspace is nongeneric at every admissible level with
denominator two.
\end{proposition}
\begin{proof}
Use the symplectic realization
\[
 x_0=\begin{pmatrix}0&0\\I_n&0\end{pmatrix},\qquad
 \n=\left\{\begin{pmatrix}A&B\\0&-A^{\mathsf T}\end{pmatrix}:
 A\text{ strictly upper triangular},\ B=B^{\mathsf T}\right\}.
\]
Then
\[
 \left[\begin{pmatrix}A&B\\0&-A^{\mathsf T}\end{pmatrix},x_0\right]
 =\begin{pmatrix}B&0\\-A^{\mathsf T}-A&-B\end{pmatrix},
\]
so $\g^{x_0}\cap\n=0$. Since $x_0^2=0$,
$(\ad x_0)^3=0$. The little adjoint representation is the primitive
part of $\Lambda^2\C^{2n}$. The operator induced by $x_0$ on the
tensor square is $x_0\otimes1+1\otimes x_0$, whose cube is zero.
Hence \eqref{eq:orbit-description} puts $x_0$ in
$\overline{\mathbb O_q}$ for every odd $q\ge3$ and every even
$q\ge4$. This proves the vanishing assertion.

For $q=2$, the orbit is $\mathbb O_{\min}$ with partition
$(2,1^{2n-2})$; see \cite{AvEM}. It is equivariantly isomorphic to
$(\C^{2n}\setminus\{0\})/\{\pm1\}$ via
$v\mapsto(w\mapsto\omega(v,w)v)$. On the standard representation,
the tangent space to the $N_+$-orbit of $(0,e_1)$ consists of
$(Be_1,-A^{\mathsf T}e_1)$ and has dimension $n+(n-1)$.
The first coordinate of the second block is an invariant linear
function, so no generic orbit has larger dimension. The finite
quotient preserves orbit dimensions, giving
$d_2=n^2-(2n-1)=(n-1)^2$.

Finally, \eqref{eq:admissibility} shows that the denominator-two
levels are exactly $k=-\frac12+j$, $j\ge0$. With suitable root-vector
normalization, the nonzero vector
\[
 \bigl(e_{2\varepsilon_1}(-1)e_{2\varepsilon_2}(-1)
       -e_{\varepsilon_1+\varepsilon_2}(-1)^2\bigr)^{j+1}\vac
\]
vanishes in the simple quotient; this is the singular-vector
construction recalled in \cite[Theorem~4.2(3)]{AMJoseph}.
It belongs to the universal principal subspace and is nonzero there
by PBW. Alternatively, the quadratic relation at $j=0$ and
Lemma~\ref{lem:propagation} suffice.
\end{proof}

\subsection{The orthogonal types}
We use the following consequence of Panyushev's canonical generic
stabilizer theory \cite[Section~2.1 and Theorem~5.3]{Pan}. For a
nilpotent orbit $\mathbb O\subset\g$, there is a reductive
algebra $\mathfrak s_\star$ such that
\begin{equation}\label{eq:panyushev-d}
 d(\mathbb O)=\frac12\bigl(\dim\mathfrak s_\star
                              -\rank\mathfrak s_\star\bigr).
\end{equation}
For $N\ge7$, the required algebras are
\begin{equation}\label{eq:panyushev-table}
\begin{array}{c|c}
\text{condition}&\mathfrak s_\star\\\hline
x^2\ne0,\ \rank x\ge\lfloor N/2\rfloor&0\\
x^2\ne0,\ r=\rank x<\lfloor N/2\rfloor&\soalg_{N-2r}\\
\text{partition }(2^{2m},1^{N-4m})
 &\soalg_{N-4m}\oplus\slalg_2^{\oplus m}.
\end{array}
\end{equation}
The generic $N_+$-stabilizer is a maximal unipotent subgroup of the
corresponding canonical reductive stabilizer, which explains
\eqref{eq:panyushev-d}. We set $\soalg_1=0$ and note that
$\soalg_2$ is a one-dimensional torus.

\begin{proposition}\label{prop:typeB}
For $\g=B_\ell$, $\ell\ge2$, one has
\[
 d_2=\left\lceil\frac\ell2\right\rceil,\qquad
 d_q=0\quad(q\ge3).
\]
The principal subspace is nongeneric at every denominator-two
admissible level $k=-\ell+\frac32+j$, $j\ge0$.
\end{proposition}
\begin{proof}
The case $\ell=2$ follows from $B_2\cong C_2$, so suppose
$\ell\ge3$. The little adjoint representation is the standard
representation of dimension $2\ell+1$. Thus
$\overline{\mathbb O_2}=\{x:x^2=0\}$, and its dense orbit has
partition
\[
 \lambda_2=\begin{cases}(2^\ell,1),&\ell\text{ even},\\
 (2^{\ell-1},1^3),&\ell\text{ odd}.
 \end{cases}
\]
Equations \eqref{eq:panyushev-d}--\eqref{eq:panyushev-table} give
$\mathfrak s_\star=\slalg_2^{\oplus m}$ if $\ell=2m$, and
$\mathfrak s_\star=\soalg_3\oplus\slalg_2^{\oplus m}$ if
$\ell=2m+1$. Hence $d_2=m$ or $m+1$, respectively.

For the remaining denominators choose $x$ with partition
\[
 \lambda_*=\begin{cases}(3,2^{\ell-2},1^2),&\ell\text{ even},\\
 (3,2^{\ell-1}),&\ell\text{ odd}.
 \end{cases}
\]
Its square is nonzero and its rank is at least $\ell$; hence
$d(Gx)=0$. Since $x^3=0$, we have $(\ad x)^5=0$ on matrices.
For odd $q\ge3$ this implies $(\ad x)^{2q}=0$, and for even
$q\ge4$ we have $x^q=0$. Thus $x\in\overline{\mathbb O_q}$
in both cases, proving $d_q=0$ by semicontinuity.

The denominator-two levels follow from \eqref{eq:admissibility}.
In Per\v se's normalization, put
\[
 Q_B=-\frac14e_{\varepsilon_1}(-1)^2
 +\sum_{i=2}^{\ell}e_{\varepsilon_1-\varepsilon_i}(-1)
                         e_{\varepsilon_1+\varepsilon_i}(-1).
\]
By \cite[Lemma~3.2 and Theorem~3.3]{PerseB}, the vector
$Q_B^{j+1}\vac$ is singular at level $-\ell+\frac32+j$ and
vanishes in the simple quotient. It is a nonzero positive-root
PBW vector. This proves nongenericity, also obtainable from $j=0$
by Lemma~\ref{lem:propagation}.
\end{proof}

\begin{proposition}\label{prop:typeD}
For $\g=D_\ell$, $\ell\ge4$, one has $d_q=0$ for every $q\ge2$.
\end{proposition}
\begin{proof}
The admissible orbit formulas \cite{ArAssoc,AvEM} give
\[
 \mathbb O_2=\mathbb O_{\lambda_2},\qquad
 \lambda_2=\begin{cases}
 (3,2^{\ell-2},1),&\ell\text{ even},\\
 (3,2^{\ell-3},1^3),&\ell\text{ odd}.
 \end{cases}
\]
In the first case the rank is $\ell$, so
$\mathfrak s_\star=0$. In the second case the rank is $\ell-1$,
so $\mathfrak s_\star=\soalg_2$. Its unipotent part is zero,
and \eqref{eq:panyushev-d} gives $d_2=0$ in both cases.
Finally, type $D$ is simply laced, and
\[
 \overline{\mathbb O_2}=\{x:(\ad x)^4=0\}
 \subseteq\{x:(\ad x)^{2q}=0\}=\overline{\mathbb O_q}
 \qquad(q\ge2).
\]
Apply Lemma~\ref{lem:semicontinuity}.
\end{proof}

\subsection{Exceptional types}
We use the orbit tables in \cite{AvEM} and Panyushev's stabilizers
\cite[Section~2.1 and Theorem~6.1]{Pan}. Formula
\eqref{eq:panyushev-d} applies here as well; in particular, a toral
canonical stabilizer has trivial unipotent part.

\subsubsection{Type \texorpdfstring{$G_2$}{G2}}\label{subsec:g2}

Let $\alpha$ be the short simple root and $\beta$ the long simple
root, and put $\gamma=2\alpha+\beta$ and $\theta=3\alpha+2\beta$.
The orbit list in \cite[Tables~1--2]{AvEM} is
\begin{equation}\label{eq:g2-orbit-table}
\begin{array}{c|c|c}
q&\mathbb O_q&\dim\mathbb O_q\\\hline
1&0&0\\
2&\widetilde A_1&8\\
3&A_1&6\\
4,5,6,9&G_2(a_1)&10\\
q\ge7,\ q\ne9&G_2&12.
\end{array}
\end{equation}
Here $A_1$ is the long-root minimal orbit and $\widetilde A_1$
is the short-root orbit.

\begin{proposition}\label{prop:g2-geometric}
For $G_2$, one has $d_3=1$, and $d_q=0$ for $q\ge2$, $q\ne3$.
\end{proposition}
\begin{proof}
Root strings give $\g^{f_\theta}\cap\n=\C e_\alpha$, so
$d(A_1)\le1$. The function $x\mapsto(x,e_\theta)$ is nonconstant
on $A_1$ and is $N_+$-invariant, since $e_\theta$ is central in
$\n$. A generic $N_+$-orbit on this six-dimensional variety
therefore has dimension at most five. Hence $d(A_1)=1$.

Next take $x=f_\gamma+f_\theta$. It is conjugate to $f_\gamma$,
since $[f_{\alpha+\beta},f_\gamma]$ is a nonzero multiple of
$f_\theta$ and $[f_{\alpha+\beta},f_\theta]=0$.
For the six positive roots, the supports of $[e_\rho,x]$ are
\[
\begin{array}{c|l}
\rho&\text{nonzero components}\\\hline
\alpha&\g_{-(\alpha+\beta)}\\
\beta&\g_{-(3\alpha+\beta)}\\
\alpha+\beta&\g_{-\alpha},\ \g_{-\gamma}\\
\gamma&\C h_\gamma,\ \g_{-(\alpha+\beta)}\\
3\alpha+\beta&\g_\alpha,\ \g_{-\beta}\\
\theta&\g_{\alpha+\beta},\ \C h_\theta.
\end{array}
\]
In a linear dependence, compare successively the components in
$\g_\alpha$, $\g_{\alpha+\beta}$, $\g_{-\alpha}$,
$\g_{-(3\alpha+\beta)}$, $\C h_\gamma$, and
$\g_{-(\alpha+\beta)}$. All coefficients vanish, so
$\g^x\cap\n=0$. Thus $d(\widetilde A_1)=0$.
The subregular and regular orbit closures contain
$\widetilde A_1$, so the assertion follows from
\eqref{eq:g2-orbit-table} and semicontinuity.
\end{proof}

To prove the converse, both denominator-three residue classes must
be treated. Fix the normalization of \cite{AL} and write
\[
 a=e_\alpha,\quad b=e_\beta,\quad c=e_{\alpha+\beta},\quad
 d=e_\gamma,\quad e=e_{3\alpha+\beta},\quad f=e_\theta.
\]
Thus $[a,b]=c$, $[a,c]=-2d$, $[a,d]=-3e$, $[b,e]=-f$,
and $[c,d]=-3f$. Write $x_r=x(-r)$, and set
\begin{equation}\label{eq:g2-Q}
 Q=\frac13d_1^2-e_1c_1+f_1a_1\in U(\mathfrak u).
\end{equation}
By \cite[Proposition~2.3]{AL}, $Q\vac$ is singular at
$k=-\frac53$. At $k=-\frac43$, the singular state is
\begin{equation}\label{eq:g2-cubic}
\begin{split}
 S=\Bigl(&\frac29d_1^3-d_1(e_1c_1-f_1a_1)\\
 &-3\bigl(e_1^2b_1-f_1e_1h_\beta(-1)-f_1^2f_\beta(-1)\bigr)
       +e_1f_2-f_1e_2\Bigr)\vac.
\end{split}
\end{equation}
The use of $f$ for the positive highest-root vector is distinct from
the notation $f_\rho$ for a negative-root vector.

\begin{lemma}\label{lem:g2-square}
In the universal affine vacuum module, with the above normalization,
\begin{equation}\label{eq:g2-descendant}
 2Q^2\vac=
 \bigl(d(-1)+e(-1)f_\alpha(0)+f(-1)f_{\alpha+\beta}(0)\bigr)S.
\end{equation}
Consequently, $Q^2\vac=0$ in $L_{-4/3}(G_2)$.
\end{lemma}
\begin{proof}
Apply the zero modes on the right by the adjoint action, using
$f_{\alpha+\beta}=[f_\beta,f_\alpha]$, and move negative modes
into PBW order. The Cartan and negative-root terms cancel. For an
explicit check, both sides of \eqref{eq:g2-descendant} have the
following ordered expression:
\begin{align*}
 &2a_1^2f_1^2-4a_1c_1e_1f_1
 +\frac43a_1d_1^2f_1+2c_1^2e_1^2
 -\frac43c_1d_1^2e_1+\frac29d_1^4\\
 &\hspace{25mm}-4d_1e_1f_2+4d_1f_1e_2-4e_1f_1d_2.
\end{align*}
Only brackets of zero modes with negative modes and brackets of
negative modes enter, so this identity is independent of the level.
At level $-\frac43$, $S$ belongs to the maximal ideal by
\cite[Propositions~2.3--2.4]{AL}, proving the last assertion.
\end{proof}

\begin{proposition}\label{prop:g2-nongeneric}
Every admissible denominator-three principal subspace in type $G_2$
is nongeneric.
\end{proposition}
\begin{proof}
The admissibility condition gives the disjoint union
\[
 \left\{-\frac53+j:j\in\Z_{\ge0}\right\}
 \sqcup\left\{-\frac43+j:j\in\Z_{\ge0}\right\}.
\]
The vectors $Q\vac$ and $Q^2\vac$ are nonzero in
$N(k\Lambda_0)$ by PBW. They vanish at the two initial levels,
respectively, by \cite{AL} and Lemma~\ref{lem:g2-square}.
Lemma~\ref{lem:propagation} covers every integral translate.
\end{proof}

\subsubsection{Types $E_6$, $E_7$, and $E_8$}

\begin{proposition}\label{prop:typeE}
Let $\mathfrak g$ be of type $E_6$, $E_7$, or $E_8$.
Then $d_q=0$ for every $q\ge2$. Hence $W(k\Lambda_0)$ is
generic at every nonintegrable admissible vacuum level.
\end{proposition}

\begin{proof}
Tables~5--7 of \cite{AvEM} identify the denominator-two orbits.
Combining these with \cite[Theorem~6.1 and Table~1]{Pan} gives
\[
\begin{array}{c|c|c|c}
 \mathfrak g & \mathbb O_2 & \mathfrak s_\star & d_2\\ \hline
 E_6 & 3A_1 & \mathfrak t_2 & 0\\
 E_7 & 4A_1 & 0 & 0\\
 E_8 & 4A_1 & 0 & 0.
\end{array}
\]
Here $\mathfrak t_2$ denotes a two-dimensional toral algebra.
Thus all three values of $d_2$ vanish by
\eqref{eq:panyushev-d}.
In particular, the nonzero canonical stabilizer in type $E_6$
does not obstruct genericity.

Since these Lie algebras are simply laced, the admissible orbit
closures satisfy (see \eqref{eq:orbit-description})
\[
 \overline{\mathbb O_q}
 =\{x\in\mathfrak g:(\operatorname{ad}x)^{2q}=0\}.
\]
Therefore
$
 \overline{\mathbb O_2}\subseteq\overline{\mathbb O_q}
 \qquad(q\ge2).
$
Lemma~\ref{lem:semicontinuity} yields
$0\le d_q\le d_2=0$.
The assertion about principal subspaces now follows from
Proposition~\ref{prop:genericity}.
\end{proof}

\subsubsection{Type $F_4$}

\begin{proposition}\label{prop:typeF}
For $\mathfrak g$ of type $F_4$, one has
\[
 d_2=9,\qquad d_q=0\quad(q\ge3).
\]
At a nonintegrable admissible vacuum level, the principal subspace
is nongeneric exactly when the denominator is two, equivalently
when
\[
 k=-\frac52+j,\qquad j\in\mathbb Z_{\ge0}.
\]
\end{proposition}

\begin{proof}
Tables~3--4 of \cite{AvEM} give
\[
 \mathbb O_2=A_1,\qquad
 \mathbb O_3=\widetilde A_2+A_1,\qquad
 \mathbb O_4=A_2+\widetilde A_1.
\]
By \cite[Theorem~6.1 and Table~1]{Pan}, their canonical
stabilizer algebras are, respectively,
$
 \mathfrak{sp}_6, 0, 0.
$
Equation~\eqref{eq:panyushev-d} therefore gives
$
 d_2=\frac{21-3}{2}=9, d_3=d_4=0.
$
Let $\varrho_{\theta_s}$ be the little adjoint representation.
Using the orbit-closure formulas 
\eqref{eq:orbit-description},
it follows that
\[
 \overline{\mathbb O_3}\subseteq\overline{\mathbb O_q}
 \quad(q\ge3\text{ odd}),\qquad
 \overline{\mathbb O_4}\subseteq\overline{\mathbb O_q}
 \quad(q\ge4\text{ even}).
\]
Semicontinuity proves $d_q=0$ for all $q\ge3$, and hence
genericity at these denominators.

It remains to treat denominator two.
Since $h^\vee=9$ and $h=12$, admissibility requires
\[
 k=-9+\frac p2,\qquad p\ge12,\qquad p\text{ odd}.
\]
Thus $p=13+2j$ and $k=-\frac52+j$, with $j\ge0$.

Use the positive-root vectors and normalization of \cite{PerseF},
and put
\[
 Q_F=-\frac14e_{\varepsilon_1}(-1)^2
 +\sum_{i=2}^{4}
 e_{\varepsilon_1-\varepsilon_i}(-1)
 e_{\varepsilon_1+\varepsilon_i}(-1).
\]
By \cite[Lemma~4.10 and Theorem~4.11]{PerseF},
$Q_F^{\,j+1}\mathbf1$ is a singular vector at level
$-\frac52+j$ and vanishes in the simple quotient.
Every current appearing in this expression belongs to
$\mathfrak n_+$. Moreover, its highest PBW symbol is the
$(j+1)$st power of a nonzero quadratic polynomial, so the vector
is nonzero in the universal principal subspace.
Consequently $W((-\frac52+j)\Lambda_0)$ is nongeneric for every $j\ge0$.
\end{proof}

\begin{remark}
The level $-\frac52$ is the lowest admissible half-integer level
for $F_4$. It is coprincipal admissible; it is not a boundary
principal admissible level in the terminology of \cite{KW}.
\end{remark}

\begin{proof}[Proof of Theorem~\ref{thm:classification}]
Propositions~\ref{prop:typeA}, \ref{prop:typeC}, \ref{prop:typeB},
\ref{prop:typeD}, \ref{prop:g2-geometric}, \ref{prop:typeE}, and
\ref{prop:typeF} give $d_q=0$ outside the listed exceptions.
Proposition~\ref{prop:genericity} proves genericity there.
The positive-root null vectors in types $B,C,F_4$ and
Proposition~\ref{prop:g2-nongeneric} establish every exceptional case.
The character formula follows from PBW.
\end{proof}

\begin{corollary}\label{cor:equivalence}
At every nonintegrable admissible vacuum level, for every simple $\g$,
\[
 W(k\Lambda_0)\text{ is generic}
 \quad\Longleftrightarrow\quad d_q=0
 \quad\Longleftrightarrow\quad
 \pi|_{\mathbb O_q}:\mathbb O_q\to\n^*\text{ is dominant}.
\]
In these cases $R_{W(k\Lambda_0)}\cong S(\n)$ and the principal
subspace is classically free.
\end{corollary}

\section{Associated varieties: a geometric conjecture}
\label{subsec:associated-variety-conjecture}
Retain the admissible level $k=-h^\vee+p/q$ and write $W=W(k\Lambda_0)$.
The projection fixed in Section~\ref{sec:classification} defines
$Y_q=\overline{\pi(\mathbb O_q)}\subset\n^*$.
The current inclusion induces a surjection
\[
 R_W\twoheadrightarrow
 \operatorname{im}\bigl(R_W\longrightarrow
 (R_{L_k(\g)})_{\mathrm{red}}\bigr)=\C[Y_q].
\]
Thus $Y_q\subseteq X_W\subseteq\n^*$, and \eqref{eq:moment-rank} gives
\[
 |\Delta_+|-d_q=\dim Y_q\le\dim X_W\le|\Delta_+|.
\]
\begin{conjecture}\label{conj:principal-associated-variety}
For every admissible level, the canonical map induces an isomorphism
\[
 (R_{W(k\Lambda_0)})_{\mathrm{red}}\cong\C[Y_q].
\]
In particular, $\dim X_{W(k\Lambda_0)}=|\Delta_+|-d_q$.
\end{conjecture}
The reduced associated variety would therefore depend only on $q$,
although the nonreduced $C_2$-structure may depend on the numerator.
This conjecture makes no assertion of classical freeness. It holds when
$d_q=0$ by Corollary~\ref{cor:equivalence}. For $q=1$, integrability makes
every positive-root generator nilpotent in $R_W$, so $X_W=Y_q=\{0\}$.
The remaining projection dimensions, already determined geometrically,
are
\[
\begin{array}{c|c|c|c}
\g&q&\mathbb O_q&\dim Y_q\\\hline
B_\ell&2&(2^{2m},1^{2\ell+1-4m}),\ m=\lfloor\ell/2\rfloor&
\ell^2-\lceil\ell/2\rceil\\
C_n&2&(2,1^{2n-2})=\mathbb O_{\min}&2n-1\\
F_4&2&A_1=\mathbb O_{\min}&15\\
G_2&3&A_1=\mathbb O_{\min}&5.
\end{array}
\]

\section{The symplectic rank-two principal subspace}\label{sec:main}
\subsection{Symplectic-boson conventions}\label{sec:free}
We first fix a realization valid for every symplectic rank. In this
section $\g=\spalg_{2n}$ and $W=W(-\frac12\Lambda_0)$.

Let $\B_n$ be the vertex algebra generated by even fields
$\beta_i,\gamma_i$, $1\leq i\leq n$, with singular operator products
\begin{equation}\label{eq:ope}
 \beta_i(z)\gamma_j(w)\sim\frac{\delta_{ij}}{z-w},\qquad
 \gamma_i(z)\beta_j(w)\sim-\frac{\delta_{ij}}{z-w},
\end{equation}
and all other pairwise OPEs regular. Give $\beta_i$ and $\gamma_i$
conformal weight $1/2$. The underlying Fock space is polynomial in
their negative creation modes. In particular,
\begin{equation}\label{eq:free-gr}
 R_{\B_n}=\C[b_1,\ldots,b_n,c_1,\ldots,c_n],\qquad
 \gr^F\B_n\cong\J(R_{\B_n}).
\end{equation}

The quadratic fields realize the affine symplectic algebra at level
$-1/2$ \cite{FF}. The subalgebra generated from the vacuum by all
quadratic currents is the simple vacuum algebra
$L_{-1/2}(\spalg_{2n})$. We use only its positive-root currents and
choose the following convenient rescaling:
\begin{equation}\label{eq:generators}
 E_{ij}=\normal{\beta_i\beta_j}\quad(1\leq i\leq j\leq n),
 \qquad D_{ij}=\normal{\beta_i\gamma_j}\quad(1\leq i<j\leq n).
\end{equation}
Their root weights are $\varepsilon_i+\varepsilon_j$ and
$\varepsilon_i-\varepsilon_j$, respectively. These fields strongly generate
$W$.

The normalization in \eqref{eq:generators} makes the rank-one symmetric
matrix relations have coefficient one. In the convention
$X_{2\varepsilon_i}=-\frac12\normal{\beta_i^2}$, one instead has
$E_{ii}=-2X_{2\varepsilon_i}$.

For the rank-two discussion below, take $n=2$ and set
\begin{equation}\label{eq:rank-two-fields}
 X=E_{11},\qquad Y=E_{22},\qquad Z=E_{12},\qquad U=D_{12}.
\end{equation}
The fields $X,Y,Z$ have mutually regular OPEs. The only nonzero
singular OPEs among these four generators, up to reversal, are
\[
 U(z)Y(w)\sim\frac{-2Z(w)}{z-w},\qquad
 U(z)Z(w)\sim\frac{-X(w)}{z-w}.
\]
There is an exact normally ordered relation
\begin{equation}\label{eq:quad-relation}
 \normal{XY}-\normal{ZZ}=0.
\end{equation}
No contraction occurs in this identity.

\subsection{The quadric and its jet algebras}\label{sec:jets}

Put $Q=\C[x,y,z]/(xy-z^2)$ and use coefficient variables
$x_r,y_r,z_r$, $r\geq0$. Let
\begin{equation}\label{eq:fr}
 f_r=\sum_{a+b=r}(x_a y_b-z_a z_b),\qquad
 Q_m=\frac{\C[x_0,\ldots,x_m,y_0,\ldots,y_m,z_0,\ldots,z_m]}
 {(f_0,\ldots,f_m)}.
\end{equation}


Next result is known; see \cite[Theorem 3.1]{Mourtada}.

\begin{lemma}\label{lem:jet-domain}
For every $m\geq0$, $Q_m$ is an integral complete intersection of
dimension $2(m+1)$. The transition maps $Q_m\to Q_{m+1}$ are
injective, and $\J(Q)$ is a domain.
\end{lemma}

\begin{remark}
The lemma is a direct proof for the $A_1$ surface singularity. It is
consistent with the general relation between canonical singularities
and irreducible jet schemes established in \cite{Mus}.
\end{remark}

\paragraph{Standard monomials.}

It is convenient to write $x[i]=x_{i-1}$, $y[i]=y_{i-1}$, and
$z[i]=z_{i-1}$, so that each of these variables has weight $i$.

\begin{lemma}\label{lem:groebner}
The residue classes of monomials not divisible by
\begin{equation}\label{eq:forbidden}
 z[i]^2\quad\hbox{or}\quad x[i+1]y[i],\qquad i\geq1,
\end{equation}
form a basis of $\J(Q)$.
\end{lemma}
\begin{proof}
Fix $m$ and use graded reverse lexicographic order on the variables
of $Q_m$, with descending variable order
\[
 y[m+1]>z[m+1]>x[m+1]>\cdots>y[1]>z[1]>x[1].
\]
For monomials of the same ordinary degree, compare exponents in the order
\[
 x[1],z[1],y[1],x[2],z[2],y[2],\ldots;
\]
the smaller exponent at
the first difference gives the larger monomial. Directly from
\eqref{eq:fr}, the leading monomials are
\[
 \operatorname{LM}(f_{2i-2})=z[i]^2,\qquad
 \operatorname{LM}(f_{2i-1})=x[i+1]y[i],
\]
whenever the indicated relation occurs among $f_0,\ldots,f_m$.
Indeed, maximal terms avoid the smallest possible indices; among the
remaining terms the stated variable order selects the displayed one.
These leading monomials are pairwise relatively prime. Buchberger's
product criterion therefore shows that $f_0,\ldots,f_m$ are a
Gr\"obner basis. The usual standard monomials form a basis of $Q_m$.

To pass to arcs, work in any fixed weight $d$. Variables of weight
greater than $d$ and relations of weight greater than $d$ cannot
affect this homogeneous component. Thus the finite standard-monomial
description stabilizes in every weight and gives \eqref{eq:forbidden}.
\end{proof}

\begin{corollary}\label{cor:quad-hilbert}
With $\wt(x_r)=\wt(y_r)=\wt(z_r)=r+1$,
\[
 \HS_q\J(Q)=\frac{1}{(1-q)(q;q)_\infty^2}.
\]
\end{corollary}
\begin{proof}
By Lemma~\ref{lem:jet-domain}, the $m+1$ equations in
\eqref{eq:fr} form a regular sequence of weights $2,3,\ldots,m+2$.
Thus
\[
 \HS_q Q_m=
 \frac{\prod_{r=0}^m(1-q^{r+2})}
 {\prod_{r=0}^m(1-q^{r+1})^3}.
\]
Taking the coefficientwise limit gives the result. Equivalently,
Lemma~\ref{lem:groebner} gives the product
\[
 \prod_{i\geq1}(1+q^i)\,
 \frac1{1-q}\prod_{i\geq1}
 \left(\frac1{1-q^i}+\frac1{1-q^{i+1}}-1\right).
\]
The factor $(1-q)^{-1}$ accounts for $x[1]$; each parenthesis counts
the mutually exclusive use of $y[i]$ and $x[i+1]$. Since the
parenthesis is $(1-q^{2i+1})/((1-q^i)(1-q^{i+1}))$, the product
simplifies to the same expression.
\end{proof}

\subsection{The \texorpdfstring{$C_2$}{C2}-algebra in rank two}

Set $R=\C[x,y,z,u]/(xy-z^2)$. The fields are normalized as in
\eqref{eq:rank-two-fields}.

\begin{proposition}\label{prop:c2}
Sending $x,y,z,u$ to the classes of $X,Y,Z,U$ induces an isomorphism
$R\cong R_{W}$. Moreover, the map $R_{W}\to R_{\B_2}$ is
injective.
\end{proposition}
\begin{proof}
Strong generation and \eqref{eq:quad-relation} give a surjection
$R\twoheadrightarrow R_{W}$. Composing with the map induced by
$W\hookrightarrow\B_2$ gives
\begin{equation}\label{eq:parametrization}
 \iota_0:R\longrightarrow\C[b_1,b_2,c_1,c_2],\qquad
 (x,y,z,u)\longmapsto(b_1^2,b_2^2,b_1b_2,b_1c_2).
\end{equation}
We show that this composite is injective.

The ring $R$ is a domain. Localizing at $x$ gives
$R[x^{-1}]=\C[x^{\pm1},z,u]$, with $y=z^2/x$.
After adjoining $b_1$ with $b_1^2=x$, the formulas
$b_2=z/b_1$ and $c_2=u/b_1$ invert the parametrization; $c_1$
is an additional free variable. This extension is faithfully flat
(indeed, free of rank two before adjoining $c_1$). Therefore the
localized map is injective. Since $x$ is a non-zero-divisor in $R$,
$\iota_0$ is injective. Both assertions follow from the factorization
$R\twoheadrightarrow R_{W}\to R_{\B_2}$.
\end{proof}

\subsection{Arc injectivity.}

\begin{proposition}\label{prop:arc-injection}
For every $m\geq0$, the map on $m$-jet algebras induced by
\eqref{eq:parametrization} is injective. In particular,
\[
 \J(R)\hookrightarrow\J\bigl(\C[b_1,b_2,c_1,c_2]\bigr).
\]
\end{proposition}
\begin{proof}
Write $R_m=Q_m[u_0,\ldots,u_m]$. By
Lemma~\ref{lem:jet-domain}, this ring is a domain.
All series below are truncated modulo $t^{m+1}$.
After inverting $x_0$, the equations
\[
 [t^r]\bigl(x(t)y(t)-z(t)^2\bigr)=0,
 \qquad 0\le r\le m,
\]
determine $y_0,\ldots,y_m$ recursively. Thus
\[
 A:=R_m[x_0^{-1}]
 =
 \mathbb C[x_0^{\pm1},x_1,\ldots,x_m,
           z_0,\ldots,z_m,u_0,\ldots,u_m].
\]
Introduce a new indeterminate $s$ and form
\[
 A'=A[s]/(s^2-x_0).
\]
The ring $A'$ is a free $A$-module with basis $1,s$,
so the natural map $A\longrightarrow A'$ is injective.
Also $s$ is invertible, with inverse $s x_0^{-1}$.
There is a unique truncated series
\[
 b_1(t)=s+\sum_{r=1}^{m}b_{1,r}t^r
 \in A'[t]/(t^{m+1})
\]
with constant coefficient $s$ and satisfying
$b_1(t)^2=x(t)$.
Indeed, its coefficients are determined successively by
\[
 b_{1,r}
 =
 \frac{1}{2s}
 \left(x_r-\sum_{a=1}^{r-1}b_{1,a}b_{1,r-a}\right),
 \qquad 1\le r\le m.
\]
Since its constant coefficient is invertible, $b_1(t)$
is invertible in $A'[t]/(t^{m+1})$. We may therefore set
$
 b_2(t)=\frac{z(t)}{b_1(t)},
 c_2(t)=\frac{u(t)}{b_1(t)}.
$
The quadric equation then gives
$
 b_2(t)^2=\frac{z(t)^2}{x(t)}=y(t).
$

Let
\[
 T_m=
 \mathbb C[b_{1,r},b_{2,r},c_{1,r},c_{2,r}:0\le r\le m]
\]
be the target jet algebra. Under the parametrization,
$x_0$ maps to $b_{1,0}^2$, so the localized map has target
$T_m[b_{1,0}^{-1}]$.
The preceding formulas give an isomorphism
\[
 T_m[b_{1,0}^{-1}]
 \cong A'[c_{1,0},\ldots,c_{1,m}],
\]
under which $b_{1,0}$ corresponds to $s$ and the map from
$A$ is its natural inclusion. The inverse is obtained
by taking coefficients in
\[
 x(t)=b_1(t)^2,\qquad
 z(t)=b_1(t)b_2(t),\qquad
 u(t)=b_1(t)c_2(t).
\]
Hence the localized jet map $R_m[x_0^{-1}] \rightarrow T_m[b_{1,0}^{-1}]$ is injective.
Since $R_m\longrightarrow R_m[x_0^{-1}]$ is injective
because $R_m$ is a domain and $x_0\ne0$.
It follows that the original map $R_m\longrightarrow T_m$
is injective. Finally, passing to the direct limit proves the
assertion for arc algebras.
\end{proof}

\begin{theorem}\label{thm:main}
The principal subspace $W=W(-\frac12\Lambda_0)\subset L_{-1/2}(\spalg_4)$ is
classically free, and
\begin{equation}\label{eq:gr-presentation}
 \gr^F W\cong
 \frac{\C[x_r,y_r,z_r,u_r:r\geq0]}
 {\left(\sum_{a+b=r}(x_a y_b-z_a z_b):r\geq0\right)}.
\end{equation}
Its character is
\[
 \ch W=\frac1{(1-q)(q;q)_\infty^3}.
\]
The embedding $W\hookrightarrow\B_2$ is strict for the Li
filtrations.
\end{theorem}
\begin{proof}
Proposition~\ref{prop:c2} identifies $R_{W}$ and its free-field
symbol map. Proposition~\ref{prop:arc-injection} proves injectivity
after applying $\J$. Since $\B_2$ is classically free by
\eqref{eq:free-gr}, Theorem~\ref{thm:criterion} gives classical
freeness and \eqref{eq:gr-presentation}. Strictness follows from
Remark~\ref{rem:strictness}.

The independent variables $u_r$, of weights $r+1$, contribute
$(q;q)_\infty^{-1}$. Combining this with
Corollary~\ref{cor:quad-hilbert}, and using finiteness of the Li
filtration in every conformal degree, gives the character.
\end{proof}

\begin{corollary}\label{cor:multichar}
Let $s$ and $t$ record the $\varepsilon_1$ and $\varepsilon_2$ weights.
Then
\[
 \ch W(s,t;q)=
 \frac{(s^2t^2q^2;q)_\infty}
 {(s^2q;q)_\infty(t^2q;q)_\infty(stq;q)_\infty(st^{-1}q;q)_\infty}.
\]
This is interpreted as a formal power series in $q$ with Laurent
polynomial coefficients in $s,t$.
\end{corollary}
\begin{proof}
The four families of variables in \eqref{eq:gr-presentation} have
root weights $2\varepsilon_1,2\varepsilon_2,\varepsilon_1+\varepsilon_2$,
and $\varepsilon_1-\varepsilon_2$. Every $f_r$ has root weight
$2\varepsilon_1+2\varepsilon_2$ and conformal degree $r+2$.
The finite-jet regular-sequence formula, followed by the
coefficientwise limit, gives the expression.
\end{proof}

\subsection{Higher level $sp(4)_{-\frac12+\ell}$ principal subspaces}

To extend the results from the previous section to levels $-\frac12+\ell$ we will need this known fact (cf. Theorem~A.4 in the
Eisenbud--Frenkel appendix to \cite{Mus}).

\begin{lemma} \label{sequence}
Let \(A\) be a nonnegatively graded \(\mathbb C\)-algebra with finite-dimensional graded components. Suppose \(f_0,f_1,\ldots\) are homogeneous of positive degrees tending to infinity, and every finite initial sequence is regular. Then the \(f_i\) are algebraically independent, and any homogeneous vector-space basis of \(A/(f_0,f_1,\ldots)\), lifted to \(A\), is a basis of \(A\) over \(\mathbb C[f_0,f_1,\ldots]\).
\end{lemma}

\begin{theorem}\label{thm:sp4-half-integral}
Let $\mathfrak g=\mathfrak{sp}_4$, and let
\[
 W_\ell:=W\bigl((\ell-\tfrac12)\Lambda_0\bigr),
 \qquad \ell\in\mathbb Z_{\ge0}.
\]
Choose positive-root vectors
\[
 x=e_{2\varepsilon_1},\qquad
 y=e_{2\varepsilon_2},\qquad
 z=e_{\varepsilon_1+\varepsilon_2},\qquad
 u=e_{\varepsilon_1-\varepsilon_2},
\]
normalized so that
\[
 [u,y]=2z,\qquad [u,z]=x,\qquad [u,x]=0.
\]
Set
\[
 P=\mathbb C[x,y,z,u],\qquad F=xy-z^2.
\]
Then the canonical maps induce isomorphisms
\[
 R_{W_\ell}\cong P/(F^{\ell+1}),
 \qquad
 J_\infty\bigl(P/(F^{\ell+1})\bigr)
 \cong \operatorname{gr}^{F}W_\ell.
\]
In particular, $W_\ell$ is classically free.

With all four currents assigned conformal degree one, its character is
\[
 \operatorname{ch}W_\ell(q)
 =
 \frac{1}{(1-q)(q;q)_\infty^3}
 \sum_{n_1,\ldots,n_\ell\ge0}
 \frac{
 q^{N_1^2+\cdots+N_\ell^2+N_1+\cdots+N_\ell}
 }{
 (q;q)_{n_1}\cdots(q;q)_{n_\ell}
 },
 \qquad
 N_i=n_i+\cdots+n_\ell.
\]
Equivalently,
\[
 \operatorname{ch}W_\ell(q)
 =
 \frac{1}{(1-q)(q;q)_\infty^3}
 \prod_{\substack{m\ge1\\
 m\not\equiv0,\pm1\pmod{2\ell+3}}}
 \frac{1}{1-q^m}.
\]
For $\ell=0$, the multiple sum and the additional product
are both interpreted as $1$.
\end{theorem}

\begin{proof}
The case $\ell=0$ is the classical-freeness theorem for
$W(-\tfrac12\Lambda_0)$ proved above. We assume $\ell\ge1$ and write
\[
 A_\ell=J_\infty\bigl(P/(F^{\ell+1})\bigr).
\]

\medskip
\noindent\emph{The diagonal map.}
By Adamovi\'c's singular-vector formula
\cite[Theorem~3.1]{Adamovic}, the state
\[
 \bigl(x(-1)y(-1)-z(-1)^2\bigr)^{\ell+1}\mathbf1
\]
vanishes in $W_\ell$. Hence there is a canonical surjection
\[
 A_\ell\twoheadrightarrow\operatorname{gr}^{F}W_\ell.
\]
The diagonal embeddings of
\cite[Theorem~3.4]{Adamovic}, iterated $\ell$ times and
restricted to principal subspaces, give
\[
 W_\ell\hookrightarrow
 W_0\otimes W(\Lambda_0)^{\otimes\ell}.
\]
The principal subspace $W(\Lambda_0)$ is classically free by
\cite[Corollary~6.4]{LiMilas}, using
$\mathfrak{sp}_4\cong\mathfrak{so}_5$.
Since classical freeness is preserved under tensor products,
we obtain a composite
\begin{equation} \label{aa}
 A_\ell
 \longrightarrow \operatorname{gr}^{F}W_\ell
 \longrightarrow
 A_0\otimes
 J_\infty(R_{W(\Lambda_0)})^{\otimes\ell}.
\end{equation}
We will prove that this composite is injective.

The algebra $R_{W(\Lambda_0)}$ is generated by the four
degree-one current classes and has no linear relations.
Consequently, the assignment
\[
 y\longmapsto s,\qquad x,z,u\longmapsto0
\]
defines a surjection of commutative algebras
\[
 R_{W(\Lambda_0)}
 \twoheadrightarrow\mathbb C[s]/(s^2).
\]
Indeed, every homogeneous defining relation has degree at
least two and therefore vanishes under this assignment.
Put
\[
 D_r=J_\infty\bigl(\mathbb C[v]/(v^{r+1})\bigr).
\]
Composing \eqref{aa} with the induced maps on arc algebras gives
\[
 \Phi_\ell:A_\ell\longrightarrow A_0\otimes D_1^{\otimes\ell}.
\]
In coefficient notation, its action on generating series is
\[
 \begin{aligned}
 x(t)&\longmapsto x(t),&
 z(t)&\longmapsto z(t),&
 u(t)&\longmapsto u(t),\\
 y(t)&\longmapsto y(t)+s_1(t)+\cdots+s_\ell(t),
 \end{aligned}
\]
where each $s_i$ comes from one of the $\ell$ tensor factors.

\medskip
Next we show injectivity after localization.
After inverting $x_0$, the series $x(t)$ is invertible.
The change of variables
\[
 v(t):=y(t)-\frac{z(t)^2}{x(t)}
\]
therefore identifies the equation $F(t)^{\ell+1}=0$
with $v(t)^{\ell+1}=0$. Thus
\[
 A_\ell[x_0^{-1}]\cong H\otimes D_\ell,
\]
where
\[
 H=
 \mathbb C[x_0^{\pm1},x_j\ (j\ge1),
                    z_j,u_j\ (j\ge0)].
\]
Similarly, $A_0[x_0^{-1}]\cong H$.
Under these identifications, the localization of $\Phi_\ell$
is
\[
 \operatorname{id}_H\otimes\delta_\ell:
 H\otimes D_\ell\longrightarrow H\otimes D_1^{\otimes\ell},
 \qquad
 \delta_\ell(v(t))=s_1(t)+\cdots+s_\ell(t).
\]
The map $\delta_\ell$ is injective. Indeed, the standard
presentation of the level-$r$ vacuum principal subspace
for $\widehat{\mathfrak{sl}}_2$ identifies it, as a
commutative differential algebra, with $D_r$.
The diagonal embedding of the level-$\ell$ vacuum module
into the tensor product of $\ell$ level-one vacuum modules
then restricts precisely to $\delta_\ell$;
see \cite[Theorems~3.1 and~5.1]{CLM}.
Hence $\Phi_\ell$ is injective after localization at $x_0$.

Finally, we have to show absence of $x_0$-torsion.
We show that multiplication by $x_0$ is injective on $A_\ell$.
Let
\[
 S_m=\mathbb C[x_j,y_j,z_j,u_j:0\le j\le m],
 \qquad
 F(x(t),y(t),z(t))=\sum_{j\ge0}F_jt^j.
\]
The quadric $xy-z^2=0$ is a rational double point and as we already discussed
its finite jet schemes, are reduced and irreducible complete
intersections. This is also true for their products with the jets
of the $u$-line,  In particular,
$
 S_m/(F_0,\ldots,F_m)
$
is a domain in which $x_0$ is nonzero. It follows that
$
 F_0,\ldots,F_m,x_0
$
is a regular sequence in $S_m$.

Use the positive grading
$\deg x_j=\deg y_j=\deg z_j=\deg u_j=j+1$.
Permuting the preceding homogeneous regular sequence
shows that the images of $F_0,\ldots,F_m$ form a regular
sequence in $S_m/(x_0)$.
The graded lifting argument therefore implies that both
$S_m$ and $S_m/(x_0)$ are free over
\[
 B_m=\mathbb C[F_0,\ldots,F_m].
\]
Consequently, the exact sequence of $B_m$-modules
\[
 0\longrightarrow S_m
 \xrightarrow{\ x_0\ }S_m
 \longrightarrow S_m/(x_0)\longrightarrow0
\]
remains exact after tensoring with any quotient of $B_m$.

Let $K_{\ell,m}\subset B_m$ be generated by the coefficients
of $t^0,\ldots,t^m$ in
\[
 (F_0+F_1t+\cdots+F_mt^m)^{\ell+1}.
\]
Then
\[
 S_m/K_{\ell,m}S_m
 =
 J_m\bigl(P/(F^{\ell+1})\bigr).
\]
The preceding exactness proves that $x_0$ is a
non-zero-divisor in this finite jet algebra.
Passing to the direct limit proves the same assertion
for $A_\ell$.
Thus
$
 A_\ell\hookrightarrow A_\ell[x_0^{-1}],
$
and the localized injectivity already proved implies that
$\Phi_\ell$ is injective. Since $\Phi_\ell$ factors through
the surjection
$A_\ell\twoheadrightarrow\operatorname{gr}^{F}W_\ell$,
that surjection is an isomorphism.
Taking Li degree zero gives
$R_{W_\ell}\cong P/(F^{\ell+1})$.

\medskip
Now it is easy to compute the character. By Lemma \ref{sequence}, used on the coefficients of $F(t)^{\ell+1}$, therefore gives
\[
 \operatorname{Hilb}A_\ell
 =
 \operatorname{Hilb}A_0\,
 \operatorname{Hilb}D_\ell,
\]
where, in this formula, the generator $v$ of $D_\ell$
has degree two. We have
\[
 \operatorname{Hilb}A_0
 =\frac{1}{(1-q)(q;q)_\infty^3},
\]
while the standard difference-condition basis for $D_\ell$
gives \cite{CLM2}
\[
 \operatorname{Hilb}D_\ell
 =
 \sum_{n_1,\ldots,n_\ell\ge0}
 \frac{
 q^{N_1^2+\cdots+N_\ell^2+N_1+\cdots+N_\ell}
 }{
 (q;q)_{n_1}\cdots(q;q)_{n_\ell}
 }.
\]
The product formula follows from the $i=1$
Andrews--Gordon identity.
\end{proof}

\section{Classical freeness in type \texorpdfstring{$G_2$}{G2}}\label{sec:g2-classical}

We return to the exceptional denominator-three family of
Section~\ref{subsec:g2} and prove the classical-freeness statement
announced in the introduction at its initial level $-5/3$.

Take $\g=G_2$ and retain the root-vector normalization
$a,b,c,d,e,f$ of Section~\ref{subsec:g2}, so that
\begin{equation}\label{eq:g2-case-brackets}
 [a,b]=c,\quad [a,c]=-2d,\quad [a,d]=-3e,\quad
 [b,e]=-f,\quad[c,d]=-3f.
\end{equation}

\subsection{ Case $\ell=-\frac{5}{3}$}
Set $V=L_{-5/3}(\g)$ and $W=W(-\frac53\Lambda_0)$.
The letter $f$ is the positive highest-root vector, not a negative-root
vector. All six currents have affine degree one.

\begin{theorem}\label{thm:g2-case-main}
Let
\[
 A=\C[a,b,c,d,e,f]/(Q),\qquad Q=\frac13d^2-ce+af,
\]
with all six generators of degree one. Then the natural maps give
\[
 R_W:=W/C_2(W)\cong A,
 \qquad J_\infty(A)\cong\gr^FW.
\]
In particular, $W$ is classically free and
\begin{equation}\label{eq:g2-case-character}
 \ch W(q)=\frac{(q^2;q)_\infty}{(q;q)_\infty^6}
 =\frac{1}{(1-q)(q;q)_\infty^5}.
\end{equation}
Here the jet variable $u^{(j)}$ has degree $j+1$.
\end{theorem}

We first record the two geometric steps used in the proof.

\begin{lemma}\label{lem:g2-case-projection}
Identify $\g\cong\g^*$ by an invariant nondegenerate form. The restriction
map
\[
 \pi:\overline{\mathbb O_{\min}}\longrightarrow\n^*
\]
has image dense in the hypersurface $Y=\Spec A$. Equivalently,
\[
 \ker\bigl(S(\n)\longrightarrow
 \C[\overline{\mathbb O_{\min}}]\bigr)=(Q).
\]
\end{lemma}
\begin{proof}
By \cite[Propositions~2.3 and~2.4]{AL}, the vector
\begin{equation}\label{eq:g2-case-null}
 s=\left(\frac13d(-1)^2-e(-1)c(-1)+f(-1)a(-1)\right)\vac
\end{equation}
vanishes in $V$. The products in this formula involve commuting pairs,
so its $C_2$ symbol is exactly $Q$.
Moreover,
\begin{equation}\label{eq:g2-case-associated}
 X_V=\overline{\mathbb O_{\min}},
\end{equation}
by the minimal-orbit theorem, recalled in \cite[Theorem~6.2]{AK}; see
also \cite{AMJoseph}. Hence $Q$ vanishes on the image of $\pi$.

Choose $\xi\in\mathbb O_{\min}$ to be a nonzero multiple of the negative
highest-root vector, normalized so that $\xi(f)=1$.
Then $\xi(a)=\xi(b)=\xi(c)=\xi(d)=\xi(e)=0$.
For $x\in\g$, use the tangent vector whose value on a linear function
$u\in\g$ is $\xi([u,x])$; changing the overall sign would not affect
its span. By \eqref{eq:g2-case-brackets}, the directions $x=e,d,c,b$
give independent variations of the coordinates $b,c,d,e$, respectively,
with coefficients $-1,-3,3,1$. Their other coordinates among $a,b,c,d,e,f$
are zero at $\xi$. Choose $h\in\mathfrak h$ with $\theta(h)=1$.
The direction $x=h$ gives a nonzero variation of $f$ and zero variations
of the other five coordinates. Thus $\rank d\pi_\xi\geq5$.

The quadratic form $Q$ has rank five on the variables $a,c,d,e,f$.
It is irreducible, and $Y$ has dimension five, with $b$ a free variable.
Since the image of $\pi$ lies in $Y$ and has dimension at least five,
its closure is $Y$. This proves the assertion.
\end{proof}


\begin{lemma}\label{lem:g2-case-jets}
Let $Y=Spec(A)$. Every finite jet scheme $J_mY$ is reduced and irreducible, and its
coordinate ring is a complete intersection. The arc algebra $J_\infty(A)$
is a domain. Moreover, the pullback induced by $\pi$ is injective:
\begin{equation}\label{eq:g2-case-arc-injection}
 J_\infty(A)\hookrightarrow
 \C[J_\infty\overline{\mathbb O_{\min}}]_{\red}.
\end{equation}
\end{lemma}
\begin{proof}[Proof of Lemma 6.3]
Write $Y=\mathbb A^1\times Z$, where $Z\subset\mathbb A^5$ is the rank-five quadric cone $Q=0$ in
$a,c,d,e,f$. Its only singularity is an ordinary double point, which is rational. By
\cite[Theorem~0.1 and Proposition~1.5]{Mus}, every $J_mZ$ is irreducible and reduced. It
therefore has the dimension $4(m+1)$ of the jets over $Z\setminus\{0\}$, so the $m+1$ jets of $Q$
form a regular sequence in $5(m+1)$ variables, and $J_mZ$ is a complete intersection. The same
holds for $J_mY=\mathbb A^{m+1}\times J_mZ$. The truncation maps are dominant, so
$J_\infty(A)=\varinjlim_m\mathbb C[J_mY]$ is a domain. Finally, \eqref{eq:g2-case-arc-injection}
follows from Lemma~6.2 and Lemma~\ref{lem:arc-dominance} applied to
$\pi\colon\overline{\mathbb O_{\min}}\to Y$.
\end{proof}

\begin{proof}[Proof of Theorem~\ref{thm:g2-case-main}]
The six currents strongly generate $W$, and \eqref{eq:g2-case-null} holds in $W$.
Therefore there is a natural surjection $A\twoheadrightarrow R_W$.
Its composite with the maps induced by $W\subset V$ is
\[
 A\twoheadrightarrow R_W\longrightarrow R_V
 \twoheadrightarrow (R_V)_{\red}
 =\C[\overline{\mathbb O_{\min}}].
\]
By Lemma~\ref{lem:g2-case-projection}, this composite is injective.
It follows that $A\cong R_W$.

For classical freeness, apply the canonical arc-to-Li-graded map and
functoriality of the Li filtration to obtain
\begin{equation}\label{eq:g2-case-chain}
 J_\infty(A)\twoheadrightarrow\gr^FW
 \longrightarrow\gr^FV
 \twoheadrightarrow(\gr^FV)_{\red}.
\end{equation}
The ambient algebra $V$ is quasi-lisse, since its associated variety
is the closure of a minimal nilpotent orbit and has only two symplectic
leaves. By \cite[Theorem~9.2]{AMArc}, its reduced singular support is
the reduced arc space of its associated variety. In terms of the
canonical maps on coordinate rings, this identifies
\[
 (\gr^FV)_{\red}
 \cong \C[J_\infty\overline{\mathbb O_{\min}}]_{\red}.
\]
Under this identification, the composite in \eqref{eq:g2-case-chain} is exactly
the pullback in \eqref{eq:g2-case-arc-injection}: both send each jet generator
to the corresponding derivative of the positive-root coordinate.
The composite is therefore injective. Since the first arrow in
\eqref{eq:g2-case-chain} is surjective, it is an isomorphism.

Finally, in $J_m(A)$ there are six variables of weight $j+1$ for each
$0\leq j\leq m$, and the coefficient of $t^r$ in $Q(u(t))$ has
weight $r+2$. Lemma~\ref{lem:g2-case-jets} gives
\[
 \operatorname{Hilb}J_m(A)
 =\frac{\prod_{r=0}^{m}(1-q^{r+2})}
 {\prod_{j=0}^{m}(1-q^{j+1})^6}.
\]
Taking the coefficientwise limit proves \eqref{eq:g2-case-character}, since
the Li filtration is finite in each conformal degree and hence
preserves the graded dimension.
\end{proof}

\begin{remark}
The argument does not assume that the ambient algebra $V$ is
classically free or that $R_V$ is reduced. Reducedness is needed for
the candidate principal arc algebra, while the quasi-lisse theorem
supplies only the reduced singular support of the ambient algebra.
Likewise, injectivity of $R_W\to R_V$ and of
$\gr^FW\to\gr^FV$ is proved by the displayed
composites, not assumed from the inclusion $W\subset V$.
\end{remark}

\subsection{Case level $-\frac43$}

Let $V=L_{-\frac43}(\frak{g})$ and $W=W(-\frac43\Lambda_0)$. Then we know from Section 3; Lemma \ref{lem:g2-square}
that $$P/(F^2) \twoheadrightarrow R_{W},$$
where  
\[
P=\mathbb C[a,b,c,d,e,f],\qquad F=\frac13d^2-ec+fa,
\]
Using the same argument as in the proof of Theorem \ref{thm:sp4-half-integral} we get
\[
\operatorname{Hilb}J_\infty(P/(F^2))
=
\frac{1}{(1-q)(q;q)_\infty^5}
\sum_{r\ge0}\frac{q^{r(r+1)}}{(q;q)_r}
=
\frac{1}{
(1-q)(q;q)_\infty^5
(q^2;q^5)_\infty(q^3;q^5)_\infty
}.
\]
This follows using the standard formula 
\[
H_{J_\infty(P/(F^2))}
=
H_{J_\infty(P/(F))}
\cdot H_{J_\infty(\mathbb C[u]/(u^2))},
\]
where $u$ has weight two.

Unfortunately, we cannot prove that $J_\infty(P/(F^2))$ is isomorphic to $W$ so we leave it here as a 

\begin{conjecture}\label{conj:rogers-ramanujian} $W(-\frac{4}{3} \Lambda_0)$ is classically free and $W(-\frac{4}{3} \Lambda_0) \cong J_\infty(P/(F^2))$.

\end{conjecture}

\section{The symplectic rank-three case}\label{sec:sp6}

This part is algebraically and combinatorially the most delicate part of the paper. 

Throughout this section, let $\g=\spalg_6$,
$W=W(-\frac12\Lambda_0)$, and $N=N(-\frac12\Lambda_0)$.
We use the free fields of Section~\ref{sec:free}, and write
\[
 (X,Y,Z,P,Q,R,A,B,D)
 =(E_{11},E_{22},E_{12},E_{13},E_{23},E_{33},D_{13},D_{23},D_{12}).
\]
Let $\mathcal A\subset W$ be the vertex subalgebra generated by the first
eight generators, omitting $D$. We determine the $C_2$-algebra and character
of $W$. The character calculation combines a lower bound from a functional
model with a degeneration of $L_{-3/2}(\mathfrak{sl}_3)$. We then give a
coefficientwise Hilbert-series criterion for classical freeness of
$\mathcal A$ and explain its consequence for $L_{-3/2}(\mathfrak{sl}_3)$

\subsection{Relations and the candidate $C_2$-algebra}
The eight quadratic null states are
\begin{equation}\label{eq:C3-quadratic-relations}
\begin{gathered}
 :XY:-:ZZ:,\quad :XR:-:PP:,\quad :YR:-:QQ:,\quad
 :XQ:-:ZP:,\\
 :ZQ:-:YP:,\quad :ZR:-:PQ:,\quad :XB:-:ZA:,\quad :ZB:-:YA:.
\end{gathered}
\end{equation}
For the following three cubic null states, write $U_i=U_{(-i)}$ and let
every displayed word act on $\vac$:
\begin{align}
 s_1&=(Q_1A_1^2-P_1A_1B_1+2Z_1A_2-2X_1B_2)\vac,
 \label{eq:C3-cubic-one}\\
 s_2&=(Q_1A_1B_1-P_1B_1^2+2Y_1A_2-2Z_1B_2)\vac,
 \label{eq:C3-cubic-two}\\
 s_3&=(Q_1R_1A_1-P_1R_1B_1+Q_1P_2-P_1Q_2)\vac.
 \label{eq:C3-cubic-three}
\end{align}
Their vanishing follows by expansion in the free fields. The derivative corrections
are part of the relations and cannot be dropped.

Let $x,y,z,p,\mathsf q,r,a,b,d$ denote the corresponding $C_2$-classes;
the symbol $\mathsf q$ distinguishes the class of $Q$ from the character
variable $q$. Put $h=pb-\mathsf q a$ and
\[
\begin{split}
 I_2=(\,&xy-z^2,\ xr-p^2,\ yr-\mathsf q^2,\ x\mathsf q-zp,\
 z\mathsf q-yp,\\
 &zr-p\mathsf q,\ xb-za,\ zb-ya\,).
\end{split}
\]
The null states give the canonical surjections
\begin{equation}\label{eq:C3-candidate-surjection}
 R_{\mathrm{cand}}=
 \frac{\C[x,y,z,p,\mathsf q,r,a,b,d]}{I_2+(ah,bh,rh)}
 \twoheadrightarrow R_W,
\end{equation}
and
\begin{equation}\label{eq:C3-arc-upper-bound}
 \J(R_{\mathrm{cand}})\twoheadrightarrow\J(R_W)
 \twoheadrightarrow\gr^F W.
\end{equation}

\subsection{The $C_2$-algebra of $W$}

In this part we determine the $C_2$ algebra $R_W$ and its associated variety.

Let \[S_8=\mathbb{C}[x,y,z,p,\mathsf q, r, a, b],\quad I_T=I_2+(h),\quad  J=I_2+(ah,bh,rh). \]
Define \[T:=S_8/I_T,\qquad \tilde{T}:=S_8/J.\] 
Thus, we have $R_{\mathrm{ cand}}=\tilde{T}\otimes \mathbb{C}[d]. $
Let $\alpha: R_{\mathrm{cand}}\twoheadrightarrow  R_W$ be the map in (\ref{eq:C3-candidate-surjection}), and Let $\varphi: R_{W}\rightarrow R_{\B_3}$ be induced by the free-field embedding. 
\begin{proposition}\label{prop:I_T}
Let $\varphi_0: S_8\rightarrow \mathbb{C}[b_1,b_2,b_3,c_3]\subset R_{\B_3}$ be the map sending 
\[
x,\quad y,\quad z,\quad p,\quad
\mathsf{q},\quad r,\quad a,\quad b
\]
to
\[
b_1^2,\quad b_2^2,\quad b_1b_2,\quad b_1b_3,\quad
b_2b_3,\quad b_3^2,\quad b_1c_3,\quad b_2c_3.
\] 
Then, $\ker\varphi_0=I_T. $
\end{proposition}
\begin{proof}
We equip $S_8$ with the graded reverse lexicographic order with $x>y>z>p>\mathsf{q}>r>a>b$. Then, the leading terms of eight quadrics generating $I_2$ and of $h$ are 
\begin{align}\label{leading terms}
xy,\quad p^2,\quad \mathsf{q}^2,\quad zp,\quad yp,\quad p\mathsf{q},\quad za,\quad ya,\quad \mathsf{q}a.
\end{align} 
Let $E$ be the monomial ideal generated by these leading terms. We compute the number of degree $n$ monomials in $S_8/E$. Put $\nu(n)=\binom{n+3}{3}$, the number of monomials of degree $n$ in four variables. We sort the degree $n$ monomials in $S_8/E$ into $4$ classes:
\begin{itemize}
\item $p$ occurs: then $y,z,\mathsf q$ cannot occur, and each monomial is $p$ times a degree $n-1$ monomial in $x,r,a,b$. The number of such monomials is $\nu(n-1)$. 
\item $p,\mathsf{q}$ do not occur, $a$ does at least once: each monomial is $a$ times a degree $n-1$ monomial in $x,a,r,b$. The number of such monomials is $\nu(n-1)$. 
\item $\mathsf{q}$ occurs, $p$ does not: then $a$ does not occur, and each monomial is $q$ times monomial in $x,y,z,r,b$ not divisible by $xy$. The number of such monomials is $\nu(n-1)+\nu(n-2)$.
\item None of $p,q,a$ occurs: then each monomial is a monomial in $x,y,z,r,b$ not divisible by $xy$. The number of such monomials is $\nu(n)+\nu(n-1).  $
\end{itemize}
 The total number of degree $n$ monomials is $\nu(n)+4\nu(n-1)+\nu(n-2)=(n+1)^3$. 
	
Denote the set of degree $2n$ monomials in $\operatorname{im} \varphi_0$ by $\mathsf S_{2n}$. Then we claim
\begin{align}\label{prop: characterize I_T}
\mathsf S_{2n}=\left\{
b_1^{e_1} b_2^{e_2} b_3^{e_3} c_3^k
\;\middle|\;
e_1+e_2+e_3+k=2n
\quad\text{and}\quad
k \leq e_1+e_2
\right\}.
\end{align}
Indeed, $\mathsf S_{2n}$ is contained in the RHS of (\ref{prop: characterize I_T}). Conversely, given any $b_1^{e_1}b_2^{e_2}b_3^{e_3}c_3^k$ on the RHS of (\ref{prop: characterize I_T}), we can pair each $c_3$ with a factor $b_1$ or $b_2$, which is possible since $k\leq e_1+e_2$, and pair the remaining factors $b_j$, whose number is even, arbitrarily. Then, we compute $|\mathsf{S}_{2n}|$, the cardinality of  $\mathsf{S}_{2n}$. Let $s=e_1+e_2$. There are $s+1$ choices of $b_1^{e_1}b_2^{e_2}$. Then, there are $\min\{s,2n-s\}+1$ choices  of $c_3^k$. So, 
\begin{align*}
|\mathsf{S}_{2n}|&=\sum_{s=0}^{2n}(s+1)(\min\{s,2n-s\}+1)\\
&=(n+1)^3.
\end{align*}

Since $I_T\subset \ker \varphi_0$ and $E\subset \operatorname{in}(I_T)$, where $ \operatorname{in}(I_T)$ is the leading-term ideal with respect to the order on $S_8$, 
\[(n+1)^3=|\mathsf{S}_{2n}|=\dim(S_8/\ker \varphi_0)_{n}\leq \dim(S_8/I_T)_n\leq \dim(S_8/E)_n=(n+1)^3.\] Consequently, $I_T=\ker \varphi_0$ and $\operatorname{in}(I_T)=E$, which implies that the nine quadrics generating $I_T$ are a Gr\"{o}bner basis.  \end{proof}
\begin{corollary}\label{cor:HilbT}
$T$ is a domain, and 
\begin{align*}
\operatorname{Hilb}T=\frac{1+4t+t^2}{(1-t)^4}.	
\end{align*}
In particular, $\dim (\operatorname{Spec}(T))=4$. 
\end{corollary}
\begin{proof}
$T$ is a subring of a polynomial ring, hence a domain. By the proof of Proposition \ref{prop:I_T}, $\dim(T_n)=\dim |\mathsf S_{2n}|=(n+1)^3$. So we easily obtain
\begin{align}\label{formu:HilbT}
\begin{split}
\operatorname{Hilb}(T)&=\sum_{n\geq 0} (n+1)^3 t^n=\frac{1+4t+t^2}{(1-t)^4}.
\end{split}
\end{align}
Since the pole of the Hilbert series at $t=1$ is $4$, the last assertion follows from the Hilbert-Serre theorem.
\end{proof}

Let $\psi: \tilde{T}\rightarrow T$ denote the natural map. 
\begin{lemma}\label{lem:kerpsi}
We have $\ker \psi=\mathbb{C}h$.
\end{lemma}
\begin{proof}
A direct calculation gives
\begin{align*}
	&p(xb-za)-a(x\mathsf q-z p)=xh,\quad r(xb-za)+a(zr-p\mathsf q)-b(xr-p^2)=ph,\\
	&p(zb-ya)-a(z\mathsf q-yp)=zh,\quad r(zb-ya)+a(yr-\mathsf{q}^2)-b(zr-p\mathsf{q})=\mathsf{q}h,\\
	&\mathsf{q}(zb-ya)-b(z\mathsf{q}-yp)=yh.
\end{align*}
It follows that $(x,y,z,p,\mathsf{q},r,a,b)h\subset J$. We also observe that $h\notin J$, and \[I_T=J+\mathbb{C}h+(x,y,z,p,\mathsf{q},r,a,b)h.\] Therefore, $I_T/J=\mathbb{C}h$. The result now follows. 
\end{proof}
By Proposition \ref{prop:I_T},  we may identify $T$ with its image in $\mathbb{C}[b_1,b_2,b_3,c_3]$. Using this realization of $T$, the composition 
\[\varphi\circ \alpha: R_{\mathrm{cand}}\twoheadrightarrow R_W \rightarrow R_{\B_3},\] can be identified with the map \[\tilde{T}[d]\longrightarrow T[b_1c_2],\] by sending $\tilde{T}$ to $T$ via $\psi$ and $d$ to $b_1c_2$.  According to Lemma \ref{lem:kerpsi}, we further have
\begin{align}\label{formula:ker(al)} \ker\alpha \subset \ker(\varphi\circ \alpha)=h\mathbb{C}[d].
\end{align}

It therefore remains to show that $hd^m\neq 0$ in $R_W$. To this end, we introduce some notations. We assign to each monomial in $R_{\mathrm{cand}}$ a bidegree whose first component is its polynomial degree and whose second component is its $\mathfrak{g}$-weight, with weights extended additively under multiplication. The weights of the generators are listed in the following table. 
\begin{table}[ht]
\centering\small\setlength{\tabcolsep}{4pt}
\label{tab:classes}
\begin{tabular}{l*{9}{c}}
\toprule
class in $R_{\mathrm{cand}}$ & $x$ & $y$ & $z$ & $p$ & $\mathsf{q}$ & $r$ & $a$ & $b$ & $d$\\
\midrule
root weight & $2\varepsilon_1$ & $2\varepsilon_2$ & $\varepsilon_1{+}\varepsilon_2$ & $\varepsilon_1{+}\varepsilon_3$
            & $\varepsilon_2{+}\varepsilon_3$ & $2\varepsilon_3$ & $\varepsilon_1{-}\varepsilon_3$
            & $\varepsilon_2{-}\varepsilon_3$ & $\varepsilon_1{-}\varepsilon_2$\\
\bottomrule
\end{tabular}
\end{table}
In particular, the bidegree of $hd^m$ is $(m+2,(m+1)\varepsilon_1+(1-m)\varepsilon_2)$. These degrees are pairwise different. Recall that $J_\infty(R_{\B_3})=\mathbb{C}[b_i^{(k)}, c_i^{(k)}]$. We define the jet order of a monomial to be the sum of $k$. For $v\in J_\infty(R_{\B_3})$, let $\sigma_j(v)$ be its component of jet order $j$. 
\begin{lemma}\label{lem:freefield}
Let $U^1,..., U^k$, with classes $u^1,...,u^k$ in $R_W$. Let $n_i\geq 1$, $\mathsf U=U^1(-n_1)\cdots U^k(-n_k)\mathbf 1$. Let $\iota$ be the embedding from $W$ to $\B_3$. Then $\overline{\iota(\mathsf{U})}$ in $\gr^F(\B_3)$ has no component of jet order less than $\sum_i(n_i-1)$.
\end{lemma}
\begin{proof}
Let $\epsilon$ and $\eta$ be two fields. Then, $U(-n)=\sum_{r+s=-n}\epsilon(r)\eta(s)$. If $\epsilon$ and $\eta$ are not conjugate, then $\epsilon(r)$ and $\eta(s)$ commute. Hence, the creation-only part of $U(-n)$ acts by multiplication by  $\frac{T^{n-1}(\epsilon \eta)}{(n-1)!}$. In every other term, exactly one factor, say $\eta(\frac12 + s)$ with $s\geq 0$, is an annihilation mode. It removes the term $\epsilon(-\frac12-s)$ by decreasing the jet order by $s$. In the expression of $U(-n)$, $\eta(\frac12+s)$ always occurs together with $\epsilon(-n-\frac12-s)$; $\epsilon(-n-\frac12-s)$ increases the jet order by $n+s$. Together, they increase the jet order by $n$. The result follows by applying the factors from right to left. 
\end{proof}

We now prove the following theorem.
\begin{theorem}\label{thm:c2ofW}
The $C_2$-algebra of $W$ is isomorphic to $R_{\mathrm{cand}}$. 	
\end{theorem}
\begin{proof}
By the preceding discussion, it suffices to show that $hd^m\neq 0$ in $R_W$. The class $hd^m$ is represented by 
$v_m:=P(-1)B(-1)D(-1)^m\mathbf 1-Q(-1)A(-1)D(-1)^m\mathbf 1\in W$, and we must show that $v_m\notin C_2(W)$. The bidegree of $hd^m$ is $(m+2,(m+1)\varepsilon_1+(1-m)\varepsilon_2$). By a direct calculation, the $C_2(W)$ is spanned, in bidegree $ (i+2,(m+1)\varepsilon_1+(1-m)\varepsilon_2)$ ($i\leq m$), by the vectors with factors $ZD^m$ or $PAD^{m-1}$ with one mode $-2$, and those with factors $XD^{m-1}$.  

We now prove the theorem by contradiction. Suppose $v_m\in C_2(W)$. By the free-field realization of $W$, it follows that $\iota(v_m)$ has jet order at least $1$. By Lemma \ref{lem:freefield}, $\sigma_1(\iota(v_m))$ lies in the span of 
\begin{align*}
&\text{from $ZD^m$:} && s_1=T(b_1b_2)(b_1c_2)^m, \quad s_2=b_1b_2\,T(b_1c_2)(b_1c_2)^{m-1},\\
&\text{from $PAD^{m-1}$:} && s_3=T(b_1b_3)\,b_1c_3\,(b_1c_2)^{m-1}, \quad s_4=b_1b_3\,T(b_1c_3) (b_1c_2)^{m-1},\\
& && s_5=b_1b_3\,b_1c_3\,T(b_1c_2)(b_1c_2)^{m-2}.   
\end{align*}
   	
Note that $\B_3\cong J_\infty(R_{\B_3})$ as vector spaces. Let $\pi$ be the ring endomorphism of $J_\infty(R_{\B_3})$ that kills every $b_3^{(k)}$ and $c_3^{(k)}$. Direct calculation gives
\begin{align*}
\pi\circ \sigma_1(\iota(v_m))&=\left[\beta_3(\frac12),\gamma_3(-\frac12)\right]\left(\beta_1(-\frac32)\beta_2(-\frac12)-\beta_2(-\frac32)\beta_1(-\frac12)\right)\left(\beta_1(-\frac12)\gamma_2(-\frac12)\right)^m\mathbf 1\\
&=\left(\beta_1(-\frac32)\beta_2(-\frac12)-\beta_2(-\frac32)\beta_1(-\frac12)\right)\left(\beta_1(-\frac12)\gamma_2(-\frac12)\right)^m\mathbf 1.
\end{align*}
We can view the right-hand side  as $M:=(Tb_1b_2-b_1Tb_2)(b_1c_2)^m\in J_\infty(R_{\B_3})$. We observe that $M$ cannot be written as a linear combination of $s_1$ and $s_2$. This is a contradiction because $\pi\sigma_1(\iota(v_m))$ should lie in the span of the $\pi(s_i)$ $(1\leq i\leq 5)$, and $\pi(s_3)=\pi(s_4)=\pi(s_5)=0$. Hence, $v_m\notin C_2(W)$, completing the proof.     	
\end{proof}
\begin{corollary}
Conjecture \ref{conj:principal-associated-variety} holds for $\g=\mathfrak{sp}_6$ at $k=-\frac12$.
\end{corollary}
\begin{proof}
By Theorem	\ref{thm:c2ofW}, $R_{W}=\tilde{T}[d].$ According to Lemma \ref{lem:kerpsi}, $\tilde{T}/(h)\cong T$. Since $T$ is a domain and $h$ is nilpotent, the nilradical of $\tilde{T}$ is $(h)$, and therefore
\begin{align*}
(R_W)_{\red}\cong T[d].	
\end{align*}
By Corollary \ref{cor:HilbT}, the associated variety has dimension $5$. 

The surjection $(R_W)_{\red}\twoheadrightarrow \C[Y_2]$ of Section \ref{subsec:associated-variety-conjecture} has prime kernel, $\C[Y_2]$ is a domain. By the table in Section \ref{subsec:associated-variety-conjecture}, $\dim Y_2=5$. A prime ideal of height zero in a finitely generated domain is zero, so the kernel vanishes. Hence the result follows. 
\end{proof}
\subsection{The character of $W$}
Here the goal is to compute the character of $W$.
We first obtain a character lower bound for $\mathcal{A}$ from a functional model. We then construct a flat degeneration of $L_{-3/2}(\mathfrak{sl}_3)$ whose special fibre surjects onto $\mathcal{A}$. The two character bounds will coincide. 
\subsubsection{A lower bound from the functional model}
Define $\widehat{S}$ to be a differential subalgebra of $J_{\infty}(R_{\B_3})$ generated by $\operatorname{im}\varphi_0$ and $Tb_1b_2-b_1Tb_2$. $\widehat{S}$ lies in the polynomial ring of the jets of the four colours $b_1,b_2,b_3,c_3$. A sector is a multidegree $\mathsf n=(k_1,k_2,k_3,l_3)$ in colours, and $\widehat{S}=\oplus_{\mathsf n} \widehat{S}_{\mathsf n}.$ Every generator has total degree $2$, so $\widehat{S}_{\mathsf n}=0$ unless $|\mathsf n|=2m$, then $\widehat{S}_{\mathsf n}$ has conformal weight $m+$ jet order. With $H_{\mathsf n}(q)=\sum_d\dim(\widehat{S}_{\mathsf n})_d q^d$, 
\begin{align*}
\operatorname{Hilb}\widehat{S}=\sum_{\mathsf n}q^m H_{\mathsf n}(q).	
\end{align*}

Recall $\mathcal{A}$ is generated by $X,Y,Z,P,Q,R,A,B$, omitting D. Let $\gr^F(\mathcal{A})$ be the associated graded algebra of $\mathcal{A}$ for the induced filtration $\mathcal{A}\cap F^p(\B_3)$, where we identify $\mathcal{A}$ with its free field image in $\B_3$. It is not hard to see that $\widehat{S}\subset \gr^F(\mathcal{A})$, and hence
\begin{align}\label{formula:lowerboundA}
	\ch\mathcal{A}=\ch\gr^F\mathcal{A} \geq \operatorname{Hilb}\widehat{S}.
\end{align}

We use a functional model as in \cite{FS} to compute $\widehat{S}_{\mathsf n}$. Fix a sector $\mathsf n=(k_1,k_2,k_3,l_3)$ with $|\mathsf{n}|=2m$. Take one variable for each factor:  
\[
\xi_1,\dots,\xi_{k_1}\ (\text{for }b_1),\qquad\eta_1,\dots,\eta_{k_2}\ (b_2),\qquad\zeta_1,\dots,\zeta_{k_3}\ (b_3),
\qquad\vartheta_1,\dots,\vartheta_{l_3}\ (c_3).
\]
Let $\Phi$ be the space of polynomials in them that are symmetric within each of the four groups. Let the dual basis of a jet monomial \[b_1^{(r_1)}\cdots b_1^{(r_{k_1})}b_2^{(s_1)}\cdots b_2^{(s_{k_2})}b_3^{(t_1)}\cdots b_3^{(t_{k_3})}c_3^{(u_1)}\cdots c_3^{(u_{l_3})} \] be the class of 
\[ \xi_1^{r_1}\cdots \xi_{k_1}^{r_{k_1}}\eta_1^{s_1}\cdots \eta_{k_2}^{s_{k_2}}\zeta_{1}^{t_1}\cdots \zeta_{k_3}^{t_{k_3}}\vartheta_{1}^{u_1}\cdots \vartheta_{l_3}^{u_{l_3}},\]
the sum of its distinct permutations within each group. Give each variable $\chi$ its own formal variable $\omega_\chi$, with colour $\kappa_\chi$. For $F\in \Phi$, expanding the arcs gives 
\begin{align}\label{formula: functional} \langle \prod_{\chi} \kappa_\chi(\omega_{\chi}),F\rangle=F((\omega_\chi)_{\chi}),\end{align}
  the value of the $F$ at the point whose $\chi$-coordinate is $\omega_\chi$, for each variable $\chi$. 
  
  An allowed matching $\mu$ is a partition of the $2m$ variables into $m$ pairs of the eight types $\{\xi,\xi\},\{\xi,\eta\},\{\eta,\eta\},\{\xi,\zeta\},\{\eta,\zeta\},\{\zeta,\zeta\},\{\xi,\vartheta\},\{\eta,\vartheta\}$. Let $L_\mu$ be the subspace where paired variables are equal. Pairs $\pi=\{\xi_a,\eta_c\}$ are called mixed; put $D_\pi=\partial_{\xi_a}-\partial_{\eta_c}$. 
\begin{lemma}\label{lem:H_n}
For $F\in \Phi$, $\langle F, \widehat{S}_{\mathsf{n}}\rangle=0$ if and only if $(\prod_{\pi\in \Theta} D_\pi F)|_{L_\mu}=0$, for every allowed matching $\mu$ and every set $\Theta$ of mixed pairs of $\mu$. Hence, if $J\subset \Phi$ is the orthogonal complement of the $\widehat{S}_{\mathsf{n}}$, then $H_{\mathsf n}=\operatorname{Hilb}(\Phi/J)$. 
\end{lemma}
\begin{proof}
Let $g$ be the variable of $S_8$, and $g(\omega)$ be the arc of $g$. Each $g$ is mapped to $\kappa \kappa'$, so $g(\omega)=\kappa(\omega)\kappa'(\omega)$, and $h(\omega)=\frac{d}{d\omega}b_1(\omega)b_2(\omega)-b_1(\omega)\frac{d}{d\omega}b_2(\omega).$ Then each element in $\widehat{S}_{\mathsf{n}}$ is spanned by some coefficients of the series $\prod_{\nu=1}^{m}g_\nu(\omega_\nu)$, over all generators $g_1,...,g_m$ of total multidegree $\mathsf{n}$. Then we apply $D_\pi$ in the variables $\omega$, for $\pi\in \Theta$ to (\ref{formula: functional}), and then set the two $\omega$'s of the $\nu$-th pair equal to $\omega_\nu$:
\begin{align*}
\langle \prod_\nu g_\nu(\omega_\nu),F\rangle=(\prod_{\pi\in \Theta} D_\pi F)|_{L_\mu}(\omega_1,...,\omega_m). 	
\end{align*}
 Every pair $(\mu,\Theta)$ arises in this way, which proves the first claim. The second holds because the pairing is perfect in each degree. 
\end{proof}

\medskip\noindent\textit{Forgetting the two $b$-colours.}
We identify $b_1$ and $b_2$, that is, we forget which of the $k=k_1+k_2$ factors of $b_1, b_2$ is which. This turns $\widehat{S}$ into the differential algebra $\mathsf{E}$ generated by $b_1^2, b_1b_3, b_3^2, b_1c_3$. For a multidegree $\mathsf{n}'=\{k_1+k_2,k_3,l_3\}$ in $(b_1,b_3,c_3)$, let $E_{\mathsf{n}'}$ be the series of $
\mathsf{E}_{\mathsf{n}'}$ by jet order. 

\begin{proposition}\label{prop:uncolouredformula}
Let $s,u,g$ record the degrees in $b_1,b_3,c_3$. Then
\begin{align*}
	\sum_{\mathsf{n}'}E_{\mathsf{n}'}(q)s^ku^{k_3}g^{l_3}=\frac{(s^2u^2;q)_\infty}{(s^2;q)_{\infty}(u^2;q)_{\infty}(su;q)_{\infty}(sg;q)_{\infty}}.
\end{align*}
\end{proposition}
\begin{proof}
Let $R_{13}=\mathbb{C}[x,r,p,a]/(xr-p^2)$. Give the variables the weights with $q^i$ for jet order $i$:
\[x_i\mapsto s^2q^i,\quad r_i\mapsto u^2q^i,\quad p_i\mapsto suq^i,\quad a_i\mapsto sgq^i.\]
The Hilbert series of $J_{\infty}(\mathbb{C}[x,r,p,a])$ is 
\[\frac{1}{(s^2;q)_{\infty}(u^2;q)_{\infty}(su;q)_{\infty}(sg;q)_{\infty}}.\]
By renaming the index 2 as 3, Theorem \ref{thm:main} shows that $\operatorname{gr}^F(\mathsf E)\cong J_{\infty}(R_{13})$. The jets $[t^n](x(t)r(t)-p(t)^2)$ of the relation have weights $s^2u^2q^n$. The result follows from the same argument as in the proof of Corollary \ref{cor:multichar}.
\end{proof}
 \medskip\noindent\textit{Restoring the colours.}
 Fix an uncolored matching $\overline{\mu}$ and a point $o\in L_{\overline{\mu}}$. whose $m$ pairs sit at distinct positions. At $o$, the $k=k_1+k_2$ particles of colour $b_1$ or $b_2$ occupy these positions, at most two at each. We record them by the divisor
$D=\sum_p D(p)p, \deg D=k,$
where $D(p)=2$ if two such particles are paired at $p$, and $D(p)=1$ if one of them is paired with a $b_3$ or a $c_3$.

To recover an allowed matching, we choose which $k_2$ of these $k$ particles are $b_2$'s. Up to the symmetry of $F$, only their positions matter, and these form a divisor $D_2\leq D$ of degree $k_2$. Let $\Phi(\eta_1,...,\eta_{k_2})$ be a symmetric polynomial. Place the variables $\eta_j$ at the points of $D_2$, counted with multiplicity, and write $\Phi(D_2)$ for the resulting value. For a set $\Theta$ of points such that $D(p)=2$ and $D_2(p)=1$, let $\eta_{j(p)}$ be the variable placed at $p\in \Theta$, and put
\begin{align*}
	(\partial_{\Theta}\Phi)(D_2)=\left(\prod_{p\in \Theta} \partial_{\eta_j(p)}\Phi\right)(D_2).	
\end{align*}
\begin{example}
	In the sector $(2,1,0,1)$, let $\overline{\mu}$ pair two $b$'s at 0 and the third $b$ with $c_3$ at $1$. Then $D=2\cdot [0]+[1]$, and $D_2$ is $[0]$ or $[1]$. If $D_2=[0]$, then $\Theta\subset\{[0] \} $.	
\end{example}
\begin{lemma}\label{lem:Phi}
	Let $\Phi$ be symmetric in $\eta_1,...,\eta_{k_2}$ of degree at most  $k_1$ in each variable. If $(\partial_{\Theta}\Phi)(D_2)=0$ for every $D_2\leq D$ of degree $k_2$ and every set $\Theta$ of $D_2$, then $\Phi=0$.
\end{lemma}
\begin{proof}
	We prove by induction on $k_2$. Let $\Phi_p=\Phi(\eta_1,...,\eta_{k_2-1},p)$. We check that $\Phi_p$ satisfies the hypothesis of the lemma for $D-p$. Let $D_2'\leq D-p$ have degree $k_2-1$, and $\Theta'$ be the set  for $D_2'$ relative to $D-p$. Then $p\notin \Theta'$ because $(D-p)(p)\leq 1$; hence $\Theta'$ belongs to $\Theta$ for $D_2'+p$ relative to $D$. So by the induction assumption, we have that
	\begin{align*}
		( \partial_{\Theta'} \Phi_p)(D_2')=(\partial_{\Theta'}\Phi)(D_2'+p)=0.
	\end{align*}
By induction, we obtain $\Phi_p=0$. By symmetry, $\Phi$ is divisible by $\prod_j(\eta_j-p)$. 

We now consider the case when $D(p)=2$. Let $\Phi'_p=(\partial_{\eta_{k_2}}\Phi)(\eta_1,...,\eta_{k_2-1},p)$. Let $D_2'$ and $\Theta'$ be as above. If $D_2'(p)=0$, then $p$ belongs to $\Theta$ for $D_2'+p$ relative to $D$. Therefore, 
\begin{align*}
	(\partial_{\Theta'}\Phi_p')(D_2')=(\partial_{\Theta'\cup\{p\}} \Phi)(D_2'+p)=0.
\end{align*} 
If $D_2'(p)=1$, some $\eta_j$ with $j<k_2$ sits at $p$. By preceding paragraph, $\Phi$ has a factor $\eta_j-p$. This factor passes through the $\partial_{\eta_{k_2}}$, which implies that $\Phi_p'$ equals $0$. With the result from the preceding paragraph and the symmetry, $\Phi$ is divisible by $(\eta_j-p)^2$.	

In summary, $\Phi$ is divisible by $\prod_{p}(\eta_1-p)^{D(p)}$, of degree $k>k_1$ in $\eta_1$. Therefore, $\Phi$ must be zero. This completes the induction. 	
\end{proof}

\begin{lemma}\label{lem:decom of F}
Let $\Lambda$ be the space of symmetric polynomials in $\eta_1,...,\eta_{k_2}$ of degree at most $k_1$ in each variable. Every $F\in \Phi$ is uniquely $F=\sum_{\lambda}f_{\lambda}m_{\lambda}(\eta)$, where $f_\lambda\in \Phi'$ and $m_\lambda(\eta)$ is a basis of the $\Lambda$ with $\lambda$ in the $k_2\times k_1$ box. 	
\end{lemma}
\begin{proof}
By \cite[Theorem 3.1]{Andrews}, \[\operatorname{Hilb}\Lambda=\frac{(q;q)_k}{(q;q)_{k_1}(q;q)_{k_2}}.\] 
So, $\operatorname{Hilb}\Phi=\operatorname{Hilb}\Phi'\cdot \operatorname{Hilb}\Lambda$. 
It suffices to show that $\sum_\lambda f_{\lambda}m_{\lambda}(\eta)=0$ forces all $f_\lambda=0$.  
Take one point $o$ from $U=\mathbb{C}^{2m}\setminus\cup_{i<j}\{x_i=x_j\}.$ Then each $b$-particle ocupies exactly one point. For each $D_2\leq D$ of degree $k_2$, we choose $k_2$ $b$-particles to be $\eta$, and the rest is $\xi$. By assumption, $\sum_{\lambda} f_{\lambda}(o)m_{\lambda}$ vanishes at every $D_2$. By Lemma  \ref{lem:Phi}, $f_\lambda(o)$ vanishes. The result follows because $U$ is dense in $\mathbb{C}^{2m}$. 
\end{proof}
Set $\binom{k}{k_1}_q=\operatorname{Hilb}\Lambda$. Let $\Phi'$ be the polynomials which are symmetric within $\{\xi_1,...,\xi_k\}$, $\{\zeta_1,...,\zeta_{k_3}\}$ and $\{\vartheta_1,...,\vartheta_{l_3}\}$. Let $I\subset \Phi'$ such that $\langle I, \mathsf{E}_{\mathsf n'}\rangle=0$. 
\begin{proposition}\label{prop:H_nlowerbound}
For any sector $\mathsf{n}$, $H_{\mathsf{n}}(q)\geq \binom{k}{k_1}_qE_{\mathsf{n}'}(q)$ coefficientwise. 	
\end{proposition}
\begin{proof}
Let 
\[ \Phi_I=\left\{\sum_\lambda f_{\lambda}m_\lambda(\eta) \mid f_\lambda\in I\right\}.\] 
By Lemma \ref{lem:decom of F}, 
\[ \operatorname{Hilb}(\Phi/\Phi_I)=\binom{k}{k_1}_q \operatorname{Hilb}(\Phi'/I)=\binom{k}{k_1}_q {E}_{\mathsf{n}'}(q).\]
By Lemma \ref{lem:H_n}, $H_{\mathsf{n}}=\operatorname{Hilb}(\Phi/J)$. So it suffices to show $J\subset \Phi_I$. 

Let $F=\sum_{\lambda}f_{\lambda} m_{\lambda}(\eta)\in J$. Fix an uncoloured matching $\overline{\mu}$ and a point in $L_{\overline{\mu}}$ whose $m$ positions are distinct. Let $D$ be the divisor recording the uncoloured $b$-particles. Let $D_2\leq D$ be the subdivisor colouring the $b_2$. Let $o'\in L_{\mu}$ such that coordinate $\eta$'s sit at $D_2$. Let $\Theta$ be the set defined above for $D_2$ relative to $D$. Lemma \ref{lem:H_n} gives 
\[0=\left(\prod_{p\in \Theta}D_{\pi(p)}F)\right)(o')=(-1)^{|\Theta|}\sum_\lambda f_{\lambda}(o)(\partial_{\Theta}m_{\lambda})(D_2).\] 
For the second equality, we use the Leibniz rule and the fact that $\prod_{p\in \Theta}D_{\pi(p')}f_{\lambda}(o')=0$. 
By Lemma \ref{lem:Phi}, $\sum_{\lambda}f_{\lambda}(o)m_{\lambda}=0,$ so $f_{\lambda}(o)=0$. Such $o$ is dense in $L_{\overline{\mu}},$ so $f_{\lambda}\in I$. 
\end{proof}

Next we need an identity for summation over all the sectors.
Set $\sum_{k_1+k_2=k}\binom{k}{k_1}_q=G_k$.
\begin{corollary}\label{cor:lowerboundS}
$\operatorname{Hilb}\widehat{S}\geq \dfrac{1}{(q;q)_{\infty}}\displaystyle\sum\limits_{r,p\geq 0}\dfrac{q^{r+p}G_r G_{r+2p}}{(q;q)_r}.	$
\end{corollary}
\begin{proof}
By Proposition \ref{prop:H_nlowerbound}, we have
\[\operatorname{Hilb}\widehat{S}=\sum_{\mathsf{n}}q^mH_{\mathsf{n}}\geq \sum_{\mathsf n'}q^{(k+k_3+l_3)/2}G_kE_{\mathsf{n}'}=\sum_{k\geq 0}G_kq^{k/2}[s^k]E(s),\]
where $E(s)$ is the series in Proposition \ref{prop:uncolouredformula} at $u=g=q^{1/2}$.
We compute
\begin{align*}
E(s)&=\frac{(s^2q;q)_{\infty}}{(s^2;q)_{\infty}(q;q)_{\infty}(sq^{\frac12};q)_{\infty}(sq^{\frac12};q)_{\infty}}=\frac{1}{1-s^2}\frac{1}{(q;q)_{\infty}(sq^{\frac12};q)_{\infty}^2}\\
&=\frac{1}{(1-s^2)(q;q)_{\infty}}\sum_{i,j}\frac{(sq^{\frac12})^{i+j}}{(q;q)_i(q;q)_j}=\frac{1}{(1-s^2)(q;q)_{\infty}}\sum_r\frac{(sq^{\frac12})^{r}}{(q;q)_r}\sum_{i+j=r}\binom{r}{i}_q\\
&=\frac{1}{(1-s^2)(q;q)_{\infty}}\sum_r\frac{G_r}{(q;q)_r}	(sq^{\frac12})^r,
\end{align*}
where the third identity follows from the $q$-binomial theorem. Then
\[[s^k]E(s)=\frac{1}{(q;q)_{\infty}}\sum_{r+2p=k}\frac{G_r q^{\frac{r}{2}}}{(q;q)_r}.\]
Combining everything, the result follows. 
\end{proof}

\begin{lemma}\label{lem:keyidentity}
For $k\geq 0$, $\displaystyle\sum_{r\geq 0}\dfrac{q^rG_rG_{r+k}}{(q;q)_r}=\dfrac{k+1}{(1-q^{k+1})(q;q)^3_{\infty}}.$
\end{lemma}
\begin{proof}
By Mehler's formula for the Rogers--Szeg\H{o} polynomials
\cite[Chapter 3, Example 9, p. 50]{Andrews},
we obtain	
\[M(x,y;t):=\sum_{r\geq 0}h_r(x)h_r(y)\frac{t^r}{(q;q)_r}=\frac{(t^2xy;q)_{\infty}}{(t;q)_{\infty}(tx;q)_{\infty}(ty;q)_{\infty}(txy;q)_{\infty}},\]
where $h_r(x)=\sum_{j=0}^r\binom{r}{j}_qx^j$ is the Rogers--Szeg\H{o} polynomial. 
The $q$-Vandermonde identity
\cite[(3.3.10), p. 37]{Andrews}, $\binom{r+k}{l}_q=\sum_{i+j=l}\binom{r}{i}_q\binom{k}{j}_qq^{(r-i)j}$ gives
\begin{align*}
h_{r+k}(x)&=\sum_{j=0}^{r+k}\binom{r+k}{j}_qx^j=\sum_{j=0}^{k}\binom{k}{j}_qx^jq^{rj}h_{r}(xq^{-j}). 
\end{align*}
Hence
\begin{align*}
\sum_{r}h_{r}(y)h_{r+k}(x)\frac{t^r}{(q;q)_r}&=\sum_{j=0}^{k}\binom{k}{j}_qx^jM(xq^{-j},y;tq^j)\\
&=M(x,y;t)\sum_{j=0}^k\binom{k}{j}x^j\frac{(t;q)_j (ty;q)_j}{(t^2xy;q)_j},	
\end{align*}
where the last identity rests on the identity
\[(aq^j;q)_\infty=\frac{(a;q)_\infty}{(a;q)_j}.\]
Letting $x=y=1$, $t=q$, above identity becomes 
\begin{align}\label{lemproof:1}
\sum_{r\geq 0}\frac{q^rG_rG_{r+k}}{(q;q)_r}=\frac{1}{(q;q)_{\infty}^3}\sum_{j=0}^k\binom{k}{j}_q \frac{(q;q)_j}{1-q^{j+1}}.	
\end{align}

Now we compute $\displaystyle \sum_{j=0}^{k}\binom{k}{j}_q \dfrac{(q;q)_j}{1-q^{j+1}}$. It equals $\displaystyle \sum_{j=0}^{k}\dfrac{(q;q)_k}{(q;q)_{k-j}(1-q^{j+1})}$, which is the coefficient of $x^k$ in $\displaystyle \dfrac{(q;q)_k}{(x;q)_{\infty}}\sum_{j\geq 0} \dfrac{x^j}{1-q^{j+1}}$ by Euler's identity \cite[(2.2.5), p. 19]{Andrews}.
On the other hand, 
\begin{align*}
\frac{1}{(x;q)_\infty}\sum_{j\geq 0}\frac{x^j}{1-q^{j+1}}=\frac{1}{(x;q)_{\infty}}\sum_{i\geq 0} \frac{q^i}{1-xq^i}=\frac{d}{dx}\frac{1}{(x;q)_{\infty}}=\sum_{k\geq 0}\frac{(k+1)x^k}{(q;q)_{k+1}},	
\end{align*}
where the last identity again follows from Euler's identity. Combining everything, we obtain
\begin{align}\label{lemproof:2}
\sum_{j=0}^{k}\binom{k}{j}_q\frac{(q;q)_j}{1-q^{j+1}}=\frac{k+1}{1-q^{k+1}}. 	
\end{align}

The result follows from (\ref{lemproof:1}) and (\ref{lemproof:2}).  
\end{proof}

\begin{theorem}\label{thm:lowerboundS}
\begin{align*}
\sum_{r,p\geq 0}\frac{q^{r+p}G_rG_{r+2p}}{(q;q)_r}=\frac{(q^2;q^2)^8_{\infty}}{(q;q)^7_{\infty}}. 	
\end{align*}
\end{theorem}
\begin{proof}
By Lemma \ref{lem:keyidentity} with $k=2p$, the left hand side is $(q;q)_\infty^{-3}$ times
\begin{align*}
\sum_{p\geq 0}\frac{(2p+1)q^p}{1-q^{2p+1}}=\sum_{n\geq 0} \sigma(2n+1)q^n=\left(\sum_{m=0}^{\infty}q^{m(m+1)/2}\right)^4=\left(\frac{(q^2;q^2)^2_{\infty}}{(q;q)_{\infty}}\right)^4. ,	
\end{align*}
where $\sigma(2n+1)=\displaystyle\sum_{d\mid 2n+1} d$, and where we used Legendre's four triangular numbers theorem
\cite[Theorem 3]{ORW} and Gauss's identity \cite[Corollary 2.10]{Andrews}.
\end{proof}
We are ready to show the main result in this subsection. 
\begin{theorem}\label{theorem:lowerbound A}
We have, coefficientwise,
\begin{align*}
\ch\mathcal{A}\geq (-q;q)_{\infty}^8, \quad\quad \ch W\geq \frac{(-q;q)_{\infty}^8}{(q;q)_{\infty}}. 
\end{align*}
\end{theorem}
\begin{proof}
By (\ref{formula:lowerboundA}), Corollary \ref{cor:lowerboundS}, and Theorem \ref{thm:lowerboundS}, we have
\begin{align*}
\ch \mathcal{A}\geq \operatorname{Hilb}\widehat{S}\geq \frac{(q^2;q^2)_{\infty}^8}{(q;q)_{\infty}^8}=(-q;q)_\infty^8.
\end{align*}
The second inequality follows from Lemma \ref{lem:free-D-factor}. 
\end{proof}

\subsubsection{The eight-generator factor and its deformation}
\begin{lemma}\label{lem:free-D-factor}
There is a graded vector-space isomorphism
\[
 \C[d_1,d_2,\ldots]\otimes\mathcal A\xrightarrow{\sim}W,
 \qquad \prod_m d_m^{j_m}\otimes v\longmapsto
 \prod_mD(-m)^{j_m}v,\qquad\deg d_m=m.
\]
In particular, $\ch W=(q;q)_\infty^{-1}\ch\mathcal A$.
\end{lemma}
\begin{proof}
The eight current generators span an ideal in $\n$, so PBW gives
surjectivity. For injectivity, filter the Fock space by the number of
$\gamma_2$ creation operators. Vectors of $\mathcal A$ involve neither
$\gamma_1$ nor $\gamma_2$. If $b_r,g_r$ are the creation variables of
$\beta_1,\gamma_2$, the leading action of $D(-m)$ is multiplication by
$\delta_m=\sum_{r+s=m-1}b_rg_s$. Over the fraction field of the remaining
creation variables these are algebraically independent: their dependence
on the $g_r$ is triangular with diagonal $b_0$. Taking the highest
$\gamma_2$-degree in a finite relation, then descending, proves independence.
\end{proof}
The eight currents have nonzero brackets
\begin{equation}\label{eq:special-bracket}
 [A,R]=-2P,\quad[B,R]=-2Q,\quad[A,P]=-X,\quad
 [A,Q]=[B,P]=-Z,\quad[B,Q]=-Y.
\end{equation}
Their Lie algebra $\mathfrak a$ is a contraction of $\slthree$.
For matrix units $E_{ij}$, set $h=E_{11}-E_{22}$,
$H=E_{11}+E_{22}-2E_{33}$, and use the lattice basis
\begin{equation}\label{eq:matrix-basis}
\begin{aligned}
 X_t&=2t^3E_{12},&Z_t&=-t^3h,&Y_t&=-2t^3E_{21},\\
 P_t&=t^2(E_{13}-E_{32}),&Q_t&=t^2(E_{23}+E_{31}),\\
 A_t&=t(E_{13}+E_{32}),&B_t&=t(E_{23}-E_{31}),&R_t&=\tfrac23tH.
\end{aligned}
\end{equation}
Its $\C[t]$-span is closed under brackets; reduction modulo $t$ gives
\eqref{eq:special-bracket}, and the trace pairings vanish at $t=0$.
Specialization means reduction of the free lattice, not evaluation of
its matrix entries.
Put $L=L_{-3/2}(\slthree)$, let $\mathscr L$ be the $\C[t]$-vertex
subalgebra of $L\otimes\C[t,t^{-1}]$ generated by this basis and the
vacuum, and set $\mathcal B=\mathscr L/t\mathscr L$.
\begin{proposition}\label{prop:flat}
Every affine-degree component of $\mathscr L$ is finite free over
$\C[t]$, of rank equal to the corresponding dimension of $L$.
The nonzero fibers are $L$, and $\mathcal B$ has current brackets
\eqref{eq:special-bracket} and character
\[
 \ch\mathcal B=\ch L=\frac1{(q;q^2)_\infty^8}=(-q;q)_\infty^8.
\]
\end{proposition}
\begin{proof}
Polynomial current brackets give finite generation in each degree by
PBW. Each component is torsion-free as a submodule of
$L\otimes\C[t,t^{-1}]$, hence free over $\C[t]$. Inverting $t$
recovers all $\slthree$ currents and proves the rank assertion.
Evaluation at a nonzero parameter is surjective and therefore an
isomorphism degreewise. The character is the boundary admissible
formula of Kac--Wakimoto \cite[Remark~1]{KW}.
\end{proof}

\color{black}

\subsubsection{The explicit free field realization}
We recall a realization of $L$ suited to this deformation. The following
form of the maximal-ideal generator will identify the kernel. Let $\mathfrak b=\mathfrak h\oplus\mathbb C E_{13}\oplus\mathbb C E_{32}\oplus\mathbb C E_{12}$ be the
Borel subalgebra for the order $1>3>2$, whose highest root is $\theta=\epsilon_1-\epsilon_2$. \begin{lemma}[\cite{PerseA}]
In $V^{-3/2}(\mathfrak{sl}_3)$, the vector
\begin{align*}
	S=E_{12}(-1)H(-1)\mathbf{1}+\frac32(E_{13}(-1)E_{32}(-1)+E_{32}(-1)E_{13}(-1))\mathbf{1}
\end{align*}	
is the singular vector generating the maximal ideal. 
\end{lemma}
{{
Let $\Pi^{1/2}(0)=M_{c,\delta}(1)\otimes \mathbb{C}[\mathbb{Z}\frac{c}{2}]$ be the half-lattice vertex algebra with Heisenberg fields $c,\delta$ and vertex operators $e^{nc/2}$, $n\in \Z$, where 
\[\langle c,c\rangle=0,\quad \langle c,\delta\rangle=2,\quad \langle \delta,\delta\rangle=3.\]
Let $\mathsf S$ be the $\beta\gamma$-system generated by even fields $b,b^*$ with $b(z)b^*(w)\sim1/(z-w)$, and put $\mathcal{F}=\Pi^{1/2}(0)\otimes \mathsf S$. 		
\begin{theorem}[\cite{Ad16,BMR,ACG}]
Define $J_x\in\mathcal{F}$, for $x$ in the basis of $\mathfrak{sl}_3$, by
\begin{align*}
	J_{12}&=e^{c}, & J_{13}&={:}b\,e^{c/2}{:}, & J_{32}&={:}b^*e^{c/2}{:}, &
	J_h&=\delta-\tfrac32c, & J_H&=-3\,{:}bb^*{:},
\end{align*}
\vspace{-3ex}
\begin{align*}
	J_{21}&={:}\bigl(-\tfrac14\delta^2+\tfrac14\partial\delta+\tfrac14b^2b^{*2}
	-\tfrac12(\partial b)b^*+\tfrac12b\,\partial b^*\bigr)e^{-c}{:},\\
	J_{23}&={:}\bigl(-\tfrac12\delta b+\partial b-\tfrac12b^2b^*\bigr)e^{-c/2}{:},\\
	J_{31}&={:}\bigl(\tfrac12\delta b^*-\partial b^*-\tfrac12bb^{*2}\bigr)e^{-c/2}{:},
\end{align*}
where $J_{ij}=J_{E_{ij}}$, $h=E_{11}-E_{22}$ and $H=E_{11}+E_{22}-2E_{33}$. Extended linearly,
these fields satisfy
\[
J_x(z)J_y(w)\sim-\frac32\,\frac{\operatorname{tr}(xy)}{(z-w)^2}+\frac{J_{[x,y]}(w)}{z-w}
\qquad(x,y\in\mathfrak{sl}_3).
\]
The induced homomorphism $\rho\colon V^{-3/2}(\mathfrak{sl}_3)\to\mathcal{F}$,
$x(-1)\mathbf 1\mapsto J_x$, has kernel the maximal ideal $I$; thus
$L_{-3/2}(\mathfrak{sl}_3)\hookrightarrow\mathcal{F}$.
\end{theorem}

\subsubsection{Specialization of the realization}
 Put $s=t^{1/2}$. \color{black} Inside $\mathcal{F}\otimes\mathbb{C}[s,s^{-1}]$ set
\[
v^{j}=s^{3j}e^{jc/2}\ (j\in\mathbb{Z}),\qquad \tilde\delta=s^6\delta,\qquad
\varphi=s(b-b^*),\qquad \varphi^*=s^{-1}(b+b^*).
\]
Their modes satisfy
\[
[\tilde\delta(m),\tilde\delta(n)]=3s^{12}m\,\delta_{m+n,0},\qquad
[c(m),\tilde\delta(n)]=2s^{6}m\,\delta_{m+n,0},\qquad
[\varphi(r),\varphi^*(r')]=2\,\delta_{r+r',0},
\]
while $[c(m),c(n)]=[\varphi(r),\varphi(r')]=[\varphi^*(r),\varphi^*(r')]=0$. Moreover
$\tilde\delta(z)v^{j}(w)\sim s^6j\,v^j(w)/(z-w)$, $c$ commutes with $v^j$, and
$\partial v=\tfrac12{:}cv{:}$, where $v=v^1$. All structure constants lie in $\mathbb{C}[s]$,
although $\varphi^*$ contains $s^{-1}$.

\begin{lemma}\label{lem:fam}
	Let $\mathcal{F}_s\subset\mathcal{F}\otimes\mathbb{C}[s,s^{-1}]$ be the $\mathbb{C}[s]$-vertex
	subalgebra generated by $v^{\pm1},c,\tilde\delta,\varphi,\varphi^*$. Then
	\[
	\mathcal{F}_s/s\mathcal{F}_s\cong\mathcal{C}\otimes\mathcal{S}',
	\]
	where $\mathcal{C}=\mathbb{C}[J_\infty(\mathbb{C}^\times\times\mathbb{C})]$ is commutative, with
	coordinates $v$ and $v_2:=-\tilde\delta\,v^{-1}$, and $\mathcal{S}'$ is the
	$\beta\gamma$-system generated by $\varphi,\varphi^*$, with $\varphi(z)\varphi^*(w)\sim2/(z-w)$.
\end{lemma}
\begin{proof}
	By a PBW vector we mean a vector of the form
	\[
	c(-n_1)\cdots c(-n_a)\,\tilde\delta(-m_1)\cdots\tilde\delta(-m_b)\,
	\varphi(-r_1)\cdots\varphi(-r_e)\,\varphi^*(-r'_1)\cdots\varphi^*(-r'_f)\,v^{j},
	\]
	where $j\in\mathbb{Z}$, all $n_i,m_i\ge1$, all $r_i,r'_i\in\tfrac12+\mathbb{Z}_{\ge0}$, and the
	indices of each kind are non-increasing. Let $P$ be the $\mathbb{C}[s]$-span of the PBW vectors.
	
	We first show that $\mathcal{F}_s=P$. Reordering creation modes only uses the commutators
	displayed above, and their coefficients lie in $\mathbb{C}[s]$. The modes of $v^{\pm1}$ commute
	with those of $c$, $\varphi$ and $\varphi^*$, and with those of $\tilde\delta$ up to terms
	$\pm s^6v^{\pm1}(\cdot)$. On $v^j$ the lattice vertex operator gives
	\[
	v^{\pm1}(z)v^j=v^{j\pm1}\exp\Bigl(\pm\sum_{k\ge1}\frac{c(-k)z^k}{2k}\Bigr),
	\]
	with no singular factor because $\langle c,c\rangle=0$. Hence $P$ is stable under all modes of
	the generators. Since $P$ contains $\mathbf 1=v^0$ and $\mathcal{F}_s$ is spanned by iterated
	modes of the generators applied to $\mathbf 1$, we obtain $\mathcal{F}_s=P$.
	
	The PBW vectors are linearly independent over $\mathbb{C}[s]$. Indeed, after inverting $s$ each
	of them is a unit multiple of a Fock vector in the modes of $c,\delta,b-b^*,b+b^*$ acting on
	$e^{jc/2}$, and these Fock vectors form a basis of $\mathcal{F}$, because
	$(b,b^*)\mapsto(b-b^*,b+b^*)$ is invertible. Therefore $\mathcal{F}_s$ is a free
	$\mathbb{C}[s]$-module on the PBW vectors, and their classes form a basis of
	$\mathcal{F}_s/s\mathcal{F}_s$. In particular no PBW vector lies in $s\mathcal{F}_s$.
	
	It follows that the vertex algebra structure of $\mathcal{F}_s/s\mathcal{F}_s$ is obtained by
	setting $s=0$ in the commutators above. At $s=0$ every commutator vanishes except
	$[\varphi,\varphi^*]$, and so does the singular part of $\tilde\delta(z)v^j(w)$. Hence
	$\mathcal{F}_s/s\mathcal{F}_s=\mathcal{C}\otimes\mathcal{S}'$, where
	$\mathcal{C}=\mathbb{C}[v^{\pm1},\partial^nc,\partial^n\tilde\delta:n\ge0]$ is commutative.
	Finally, the relation $\partial v=\tfrac12cv$ allows the triangular substitution
	$c=2v^{-1}\partial v$, which gives
	\[
	\mathcal{C}=\mathbb{C}[v^{\pm1},\partial v,\partial^2v,\dots,\tilde\delta,\partial\tilde\delta,\dots]
	=\mathbb{C}[J_\infty(\mathbb{C}^\times\times\mathbb{C})]. \qedhere
	\]
\end{proof}

Direct calculation gives the following result. 
\begin{lemma}\label{lem:lim}
	Under $\rho$, extended by $t\mapsto s^2$ \color{black}, the lattice generators (\ref{eq:matrix-basis}) become
	\begin{gather*}
		X_t=2v^2,\qquad P_t={:}\varphi\,v{:},\qquad A_t={:}\varphi^*v{:},\qquad
		R_t=\tfrac12{:}\varphi^2{:}-\tfrac12s^4{:}\varphi^{*2}{:},\qquad Z_t=-\tilde\delta+\tfrac32s^6c,
	\end{gather*}
	\vspace{-3ex}
	\begin{align*}
		Y_t&={:}\bigl(\tfrac12\tilde\delta^2-\tfrac12s^6\partial\tilde\delta-\tfrac1{32}s^8\varphi^4
		+\tfrac1{16}s^{12}\varphi^2\varphi^{*2}
		+\tfrac12s^{12}\bigl((\partial\varphi)\varphi^*-\varphi\,\partial\varphi^*\bigr)
		-\tfrac1{32}s^{16}\varphi^{*4}\bigr)v^{-2}{:},\\
		Q_t&={:}\bigl(-\tfrac12\tilde\delta\varphi+s^6\bigl(\partial\varphi+\tfrac18\varphi^2\varphi^*\bigr)
		-\tfrac18s^{10}\varphi^{*3}\bigr)v^{-1}{:},\\
		B_t&={:}\bigl(-\tfrac12\tilde\delta\varphi^*+\tfrac18s^2\varphi^3
		+s^6\bigl(\partial\varphi^*-\tfrac18\varphi\varphi^{*2}\bigr)\bigr)v^{-1}{:}.
	\end{align*}
	In particular, they lie in $\mathcal{F}_s$. Modulo $s$ they become
	\begin{align}\label{items:8generatorof A}
	X\equiv2v^2,\quad Z\equiv vv_2,\quad Y\equiv\tfrac12v_2^2,\quad P\equiv v\varphi,\quad
	Q\equiv\tfrac12v_2\varphi,\quad A\equiv v\varphi^*,\quad B\equiv\tfrac12v_2\varphi^*,\quad
	R\equiv\tfrac12\varphi^2 .
	\end{align}
	Under
	\[
	\beta_1=\sqrt2\,v,\qquad\beta_2=v_2/\sqrt2,\qquad\beta_3=\varphi/\sqrt2,\qquad
	\gamma_3=\varphi^*/\sqrt2,
	\]
	for which $\beta_3(z)\gamma_3(w)\sim1/(z-w)$, these are exactly the currents $\mathcal{A}$.
\end{lemma}
\color{blue}}}

\subsubsection{The character formula}
Now, we can compare the special fibre with $\mathcal{A}$.
\begin{theorem}\label{theorem:character of A}
The correspondence between the eight specialized currents and their
free-field counterparts induces an isomorphism
$\mathcal B\cong\mathcal A$ of graded vertex algebras. Consequently, 
\begin{align*}
	\ch W=\frac{(-q;q)_{\infty}^8}{(q;q)_{\infty}}.
\end{align*}	
\end{theorem}
\begin{proof}
The map $\rho: L\mapsto \mathcal{F}$ can be extended to $L\otimes \mathbb{C}[t,t^{-1}]\rightarrow \mathcal{F}\otimes \mathbb{C}[s,s^{-1}]$ by sending $t$ to $s^2$. Since $\rho(\mathscr L)\subset \mathcal{F}_s$, $\rho$ descends to a homomorphism
\begin{align*}
\mathcal{B}\longrightarrow \mathcal{F}_s/s\mathcal{F}_s.
\end{align*}
 According to Lemma \ref{lem:lim}, the eight lattice currents (\ref{eq:matrix-basis}) are sent to eight reductions in (\ref{items:8generatorof A}), respectively, and the image is isomorphic to $\mathcal{A}$. Thus, we obtain a surjective homomorphism of vertex algebras $\phi: \mathcal{B}\rightarrow \mathcal{A}$. Meanwhile, according to Theorem \ref{theorem:lowerbound A} and Proposition \ref{prop:flat}, we have 
 \begin{align*}
 \ch(\mathcal{A})\geq (-q;q)_\infty^8=\ch \mathcal{B}.	
 \end{align*}
Hence $\phi: \mathcal{B}\twoheadrightarrow \mathcal{A}$ is an isomorphism. The second assertion then follows from Lemma \ref{lem:free-D-factor}. 
\end{proof}
\subsection{A criterion for classical freeness}
Finally, we discuss the classical freeness of $L$ and $W(-\frac12\Lambda_0)$.
\begin{proposition}
The subalgebra $\mathcal{A}$ of $W$ is classically free if and only if $\operatorname{Hilb}J_{\infty}(R_{\mathcal{A}})\leq (-q;q)^8_{\infty}$ coefficientwise. If this  condition holds, $L_{-\frac32}(\mathfrak{sl}_3)$ is classically free. 	
\end{proposition}
\begin{proof}
Combining the surjection
$J_{\infty}(R_{\mathcal{A}})\twoheadrightarrow
\gr^F(\mathcal{A})$
with Theorem~\ref{theorem:character of A}, we obtain
\[
\operatorname{Hilb}\bigl(J_{\infty}(R_{\mathcal{A}})\bigr)
\geq \ch(\mathcal{A})
\geq (-q;q)_{\infty}^{8}.
\]
The first assertion follows. 

Define $R_{\mathscr{L}}=\mathscr{L}/C_{2}(\mathscr{L})$, with $C_2(\mathscr{L})$ the $\C[t]$-span of the $a_{(-2)}b$ for $a,b\in L$. Then $R_{\mathscr{L}}$ is a quotient of $\C[t][x_1,...,x_8]$ by homogeneous ideal with $\deg(x_i)=1$. Fix a weight $\Delta$ and let $M=J_{\infty}(R_{\mathscr{L}})_{\Delta}$. There are only finitely many monomials in the $x^{(r)}_i$ that have weight $\Delta$, so $M$ is finitely generated over $\C[t]$. Write $M\cong \C[t]^r\oplus T$ with $T$ the torsion part. According to Proposition \ref{prop:flat}, $\dim(J_{\infty}{R_{L}})_{\Delta}=r$, and 
\[\dim J_{\infty}(R_{\mathcal{B}})_{\Delta}=\dim M/tM=r+\dim T/tT\geq r.\] Now we assume that $\mathcal{A}$ is classically free. By Theorem \ref{theorem:character of A}, so is $\mathcal{B}$. Therefore, we have 
\begin{align*}
	\ch L\leq \operatorname{Hilb}(J_{\infty}(R_L))\leq \operatorname{Hilb}(J_{\infty}(R_{\mathcal{B}}))=\ch \mathcal{B}=\ch L. 
\end{align*}
Hence $\operatorname{Hilb}(J_{\infty}(R_L))=\ch L$ gives the classical freeness for $L $. 
\end{proof}
\begin{remark}
Classical freeness of the principal subspace is a sufficient condition for the classical freeness of the $L_{-3/2}(\mathfrak{sl}_3)$, because we use the degeneration in one direction: generic fibre for $L_{-3/2}(\mathfrak{sl}_3)$, and special fibre for $\mathcal{A}$.	
\end{remark}

\color{black}

\section{Failure of classical freeness in symplectic rank at least four}
\label{sec:higher-rank}
For comparisons between ranks only, write $W^{(r)}$ and $N^{(r)}$ for
$W(-\frac12\Lambda_0)$ and $N(-\frac12\Lambda_0)$ with
$\g=\spalg_{2r}$. Let $\rho^{(r)}:N^{(r)}\to W^{(r)}$ be the
canonical quotient, and let $\psi^{(r)}$ be its principal subspace's
canonical arc map. Use the free-field currents from
\eqref{eq:generators}.

\subsection{A classical-freeness obstruction in rank four}

\begin{lemma}[Computer-assisted]
\label{lem:C4-computation}
For \(W^{(4)}\), one has
\[
\begin{array}{c|rrrrr}
N&0&1&2&3&4\\ \hline
\dim(W^{(4)})_N&1&16&119&645&2859\\
\dim J_\infty(R_{W^{(4)}})_N&1&16&119&645&2861.
\end{array}
\]
In particular,
\[
\dim(\ker\psi^{(4)})_4=2.
\]
\end{lemma}

\subsection{Passage to higher rank}

We next explain why the rank-four obstruction persists in every
higher rank. The relevant property is stronger than the existence
of a vertex-algebra embedding: the smaller principal subspace is
a face of the positive-root grading.

For type \(C_n\), take
\[
\alpha_i=\varepsilon_i-\varepsilon_{i+1}\quad(i<n),
\qquad
\alpha_n=2\varepsilon_n,
\]
and let
\[
Q_+=\sum_{i=1}^n\mathbb Z_{\geq0}\alpha_i.
\]
The algebra \(W^{(n)}\) is \(Q_+\)-graded, with
\[
\operatorname{wt}_{\mathrm{root}}E_{ij}
=\varepsilon_i+\varepsilon_j,
\qquad
\operatorname{wt}_{\mathrm{root}}D_{ij}
=\varepsilon_i-\varepsilon_j.
\]

\begin{lemma}\label{lem:terminal-Levi-face}
Let \(1\leq m\leq n\), and set
\[
F_m=
\sum_{i=n-m+1}^{n}\mathbb Z_{\geq0}\alpha_i.
\]
The sum of the root-homogeneous components of \(W^{(n)}\) with
weights in \(F_m\) is naturally isomorphic to \(W^{(m)}\).
Moreover, restriction to this face identifies
\[
(R_{W^{(n)}})_{F_m}\cong R_{W^{(m)}},
\qquad
(\gr^F W^{(n)})_{F_m}
\cong\gr^F W^{(m)},
\]
and
\[
J_\infty(R_{W^{(n)}})_{F_m}
\cong J_\infty(R_{W^{(m)}}).
\]
Under these identifications, the restriction of \(\psi^{(n)}\)
is \(\psi^{(m)}\).
\end{lemma}

\begin{proof}
A sum of elements of \(Q_+\) belongs to \(F_m\) only if every
summand belongs to \(F_m\).
Therefore a current monomial with root weight in \(F_m\) uses
only positive-root currents from the terminal \(C_m\) subsystem.
In the free-field realization these are precisely the currents
involving the last \(m\) symplectic-boson pairs.
This identifies the face subalgebra with \(W^{(m)}\).

The same positivity argument applies to every term spanning
\(C_2(W^{(n)})\) and every term defining the Li filtration:
if its total root weight lies in \(F_m\), all its factors lie
in \(F_m\). Translation preserves root weight.
Thus taking the face commutes with the \(C_2\)-quotient and with
the Li associated graded.

Finally, choose a root-homogeneous polynomial presentation of
\(R_{W^{(n)}}\). Its jet variables have the same root weights as
their underlying generators. A polynomial of root weight in
\(F_m\) can involve only face variables, and a homogeneous
relation outside the face cannot contribute to a relation
inside it. Hence taking the face also commutes with forming
the arc algebra. Naturality gives the last assertion.
\end{proof}

Using Lemma \ref{lem:C4-computation} we now obtain
\begin{theorem}
\label{thm:Cn-not-classically-free}
For every \(n\geq4\), the principal subspace
\[
W(-\frac12\Lambda_0)\subset L_{-1/2}(\spalg_{2n})
\]
is not classically free.
\end{theorem}

\begin{proof}
By Lemma~\ref{lem:C4-computation}, the map \(\psi^{(4)}\) has a
nonzero kernel in conformal weight four.
Apply Lemma~\ref{lem:terminal-Levi-face} with \(m=4\).
The restriction of \(\psi^{(n)}\) to the terminal \(C_4\) face
is exactly \(\psi^{(4)}\). Its nonzero kernel therefore remains
nonzero in \(\ker\psi^{(n)}\).
\end{proof}




Together with the rank-one and rank-two results, these findings
leave rank three as the remaining case in the proposed
classification:
\[
W(-\frac12\Lambda_0)\subset L_{-1/2}(\spalg_{2n})
\text{ is classically free}
\quad\Longrightarrow\quad n\leq3.
\]

\section{The principal subspace at level \texorpdfstring{$-\frac52$}{-5/2} in type \texorpdfstring{$F_4$}{F4}}
\label{subsec:F4-conjecture}

Set
\[
 V=L_{-5/2}(F_4),\qquad W=W(-\frac52\Lambda_0),\qquad\g=F_4.
\]
The associated variety of $V$ is the minimal nilpotent orbit closure
\cite{ArAssoc,AMJoseph}. Let
\[
 Y=\overline{\pi(\mathbb O_{\min})}\subset\mathfrak n_+^*,
 \qquad \pi(x)(u)=(x,u).
\]
Since $\dim\mathfrak n_+=24$ and
$d(\mathbb O_{\min})=9$, the moment-map rank formula gives
$\dim Y=15$.

The highest-root grading identifies
\[
 \mathfrak n_+
 =\mathfrak n_+(C_3)\oplus E\oplus\mathbb Ce_\theta,
 \qquad E\cong\Lambda^3_0\mathbb C^6.
\]
Here $E$ is a fourteen-dimensional symplectic representation.
Choose coordinates $x_1,\ldots,x_9$ on
$\mathfrak n_+(C_3)^*$, coordinates $v_1,\ldots,v_{14}$ on
$E^*$, and the highest-root coordinate $z$.
On the minimal-orbit chart $z\ne0$, the Joseph relations express
\[
 x_i=\frac{\nu_i(v)}{z}\qquad(1\le i\le9),
\]
where the $\nu_i$ are the corresponding quadratic moment-map
components; see \cite[Section~3.2]{BMR}.
Consequently,
\begin{equation}\label{eq:F4-projection-algebra}
 \C[Y]\cong P/I,\qquad
 P=\C[z,v_1,\ldots,v_{14},x_1,\ldots,x_9],\qquad
 I=(zx_i-\nu_i(v):1\le i\le9):z^\infty.
\end{equation}

Exact rational computations in polynomial degrees two and three give
\[
 \dim I_2=15,\qquad
 \dim I_3=346,\qquad
 \dim(P_1I_2)=343.
\]
Thus $I$ has fifteen minimal quadratic generators and three
minimal cubic generators. These calculations do not exclude
additional generators in higher degrees.

\begin{conjecture}\label{conj:F4-presentation}
The canonical surjection
\[
 R_W\longrightarrow\C[Y]
\]
is an isomorphism. Moreover, the ideal $I$ in
\eqref{eq:F4-projection-algebra} is generated by $I_2$ and three
representatives of a basis of $I_3/P_1I_2$.
\end{conjecture}


The shared Joseph-ideal construction motivates the stronger working
conjecture that $W$ is classically free, as in type $G_2$. The additional
relations above make the proof more involved but are not an obstruction.
This assertion remains separate from Conjecture~
\ref{conj:F4-presentation}.
If Conjecture~\ref{conj:F4-presentation} holds and
$J_\infty(\C[Y])$ is a domain, the dominant projection
$\mathbb O_{\min}\to Y$, together with the reduced
singular-support theorem \cite{AMArc}, gives injectivity of
\[
 J_\infty(\C[Y])
 \twoheadrightarrow\gr^F W
 \longrightarrow
 \C[J_\infty\overline{\mathbb O_{\min}}]_{\mathrm{red}},
\]
and hence classical freeness.
As a support to $ J_\infty(\C[Y]) \cong W(-\frac52 \Lambda_0)$, using Macaulay2 we showed that 
\[
 \operatorname{Hilb}_q J_\infty(\C[Y]) \equiv \operatorname{Hilb}_q {\rm gr}^F W \mod  q^4.
 \]

\color{black}

\end{document}